\documentstyle{amsppt}
\hoffset =.1in
\NoBlackBoxes
\magnification=1200
\nologo

\baselineskip18pt
\pageheight{24.5truecm}
\pagewidth{15truecm}
\vcorrection{-1truecm}

\font\eit=cmti8


\nopagenumbers
\topmatter
\rightheadtext{\eit Coefficients of orthogonal series}

\title On complete characterization of coefficients  of a.e. converging orthogonal series \endtitle
\author Adam Paszkiewicz (\L\'od\'z)\endauthor
\abstract We characterize sequences of numbers $(a_n)$ such that $\sum_{n\geq 1} a_n\Phi_n$ converges a.e. for any orthonormal system $(\Phi_n)$ in any $L_2$-space. In our
criterion, we use the set $B =\{\sum_{m\geq n} |a_m|^2; n\geq 1\}$ and its information function
$$h_B(t) = -\log_3(\beta-\alpha)$$
for $t\in (\alpha, \beta]$, $[\alpha, \beta]\cap B =\{\alpha, \beta\}.$

\endabstract

\endtopmatter

\footnote""{\it 2001 Mathematical Subject Classification:  Primary 40A30, 42C15. Secondary 60G07, 94A17.}

\flushpar{\bf 1. Introduction.}  The aim of this paper is to give a complete characterization of sequences $(a_n)$ for which
$$\sum a_n \Phi_n \quad \text{converges a.e. for any orthonormal sequence $(\Phi_n)$}
\tag{*}$$
(O.N. for short) in any $L_2$ space.

The celebrated Rademacher-Menshov theorem gives complete characterization of so called Weil coefficients $r_n$, $n\geq 1$, for which the convergence
$$\sum r_n |a_n|^2 <\infty, \quad a_n\in\Bbb C$$
implies (*).
Namely increasing $r_n$'s {\it are} Weil coefficients if the sequence $r_n/\log_2^2n$ is bounded and {\it are not} Weil coefficients if $r_n/\log_2^2 n\to\infty$.

The theorem, however, does not give complete characterization of sequences $(a_n)$ such that (*) is satisfied.
Some weaker sufficient conditions was  obtained by Talagrand [5], Moric and Tandori [3], and by Weber [6]. But the problem of a complete charcayerization was still open, as was stressed, in particular, in [5].

 Let us assume for simplicity that  $\sum |a_n|^2 =1$, and denote by $h_B$ the information function of the partition of $(0, 1]$ given by $B = \{\sum_{m\geq n} a_m^2; n\geq 1\}$. More precisely, $h_B: (0, 1]\to\Bbb R^+$ is given by
$$h_B(t) = -\log_3(\beta -\alpha) \quad\text{for } t\in (\alpha, \beta]$$
with $[\alpha, \beta]\cap B = \{\alpha, \beta\}$.

For  technical reasons we use $\log_3$ instead of (the more standard) $\log_2$. 
In the whole paper $|| \ ||$ denotes $L_2$-norm. Sometimes $|| \ ||$ is used on $L_0^+$-spaces, the infinite value of $|| \ ||$ is then possible.

In our investigation of information function the following notions are crucial:
\medskip

\flushpar{\bf 1.1. Notation.} In the space $(0, 1]$ with Lebesgue measure $\lambda$, let $\Cal F_i$ be the $\sigma$-field generated by
$$\{(0, n 3^{-2^i}]; \quad 1\leq n\leq 3^{2^i}\}$$
and let $|| \ ||_i$ be the conditional $L_2$-norm
$$||h||_i = (\Bbb E (h^2 |\Cal F_i))^\frac12,\quad h\geq 0,$$
taking values in $[0, \infty]$

We also put
$$V h =\lim_{i\to\infty} ||V_0\dots V_i h||$$
with
$$V_j h = (h\wedge 2^j) + ||(h -2^j)^+||_j.$$
\baselineskip-12pt
\hfill $\square$

\baselineskip18pt
\medskip

\proclaim{1.2. Theorem } For $\sum_{n\geq 1} |a_n|^2 =1$, $B =\{\sum_{m\geq n} |a_m|^2; n\geq 1\}$ and $h_B(t) =-\log_3(\beta-\alpha)$ if $t\in (\alpha, \beta]$, $[\alpha, \beta] \cap B =\{\alpha, \beta\}$, we have:

\medskip

{\rm A.} Condition (*) is equivalent to $V h_B <\infty.$

{\rm B.} The existence of an O.N.-system $(\Phi_n)$ in $L_2[0, 1]$ with
$\sum a_n \Phi_n$ diverging a.e. is equivalent to  $V h_B=\infty.$
\endproclaim

The proof requires many steps (Sections 2, 3, 4 and 5).
Applications of Theorem 1.2 are presented in Section 6. In particular, an a.s. continuity of processes with orthogonal increments in $L_2$ space of random variables is described there.

It is worth to compare our condition $V h_B <\infty$ with some formulations of Tandori criterion of unconditional convergence of orthogonal series, and with Moric Tandori criterium of convergence of orthogonal series $\sum a_n \Phi_n$ with decreasing coefficients $a_n$.

Namely, by the classical Tandori theorem, 
$$\sum a_{\sigma(n)}\Phi_n \  \text{ converges a.e. for any permutation } (a_{\sigma(n)})
\quad \text{of } \ (a_n)\tag{**}$$
 and any O.N.-system  $(\Phi_n)$
if only
$$\sum_i (\sum_{2^{2^i}\leq n < 2^{2^{i+1}}} |a_n|^2 \log_2^2 n)^{\frac12} <\infty.$$
Moreover, there exists an O.N.-system $(\Phi_n)$ such that for any sequence $\lambda_i \geq 0$, $\sum_{i\geq 1} \lambda_i =\infty$, there exist numbers $(a_n)$ and permutation $a_{\sigma(n)}$, $n\geq 1$, satisfying
$$(\sum_{2^{2^i}\leq n <2^{2^{i+1}}}|a_n|^2 \log_2^2 n)^\frac12 =\lambda_i,\quad i\geq 1,$$
and
$$\sum a_{\sigma(n)} \Phi_n \quad\text{diverges a.e.}$$
(see [2], and [4] for simplified proof).

It is known that the following conditions, formulated by the use of {\it distribution of magnitude} of $|a_n|$ only, can be obtained from  the Rademacher-Menshov and Tandori theorems in a rather simple way (see [3]).

\proclaim{1.3. Theorem } For  decreasing modules $|a_n|$,  condition (*) is equivalent to
$$\sum_{n\geq 1} |a_n|^2 \log_2^2 |a_n|<\infty.\tag{$\alpha$}$$
\endproclaim

\proclaim{1.4. Theorem } For any sequence $(a_n)$,  condition (**) is equivalent to
$$\sum_{i\geq 1} (\sum_{\underset{2^{-2^{i+1}} \leq |a_n| < 2^{-2^i}}\to{n\geq 1}} |a_n|^2 \log_2^2|a_n|)^\frac12 <\infty.\tag{$\beta$}$$

\endproclaim

Assuming, for simplicity, that $a_n\geq 0$ for $n\geq 1$, we have another formutation of $(\beta)$.

\proclaim{1.5. Proposition} For $a_n \geq 0$, condition $(\beta)$ is equivalent to
$$\sum_{i\geq 1} (\sum_{n\geq 1} a_n^2 (-\log_2 a_n)_i^2)^\frac12 <\infty \tag{$\gamma$}$$
where we use the notation 
$$z_i = z\wedge 2^{i+1} - z\wedge 2^i,\quad i\geq 1,$$
for any positive number $z.$

\endproclaim

\demo{Proof} Relations between formulas of type $(\beta)$ and $(\gamma)$ are well-known but we recall here an elementary proof.

Let $u_i = 1_{\{n; 2^i\leq -\log_2 a_n < 2^{i+1}\}} \in L_2(\Bbb N, 2^{\Bbb N}, \mu)$, on measure space $\Bbb N = \{1,2,\dots\}$ with $\mu =\sum_n a_n^2 \delta_n$. In this formula, $\delta_n$ is a  Dirac measure concentrated in $n$. Then  $u_1, u_2,\dots$ are
orthogonal vectors in $L_2$ with $\sum ||u_i||^2 <\infty$, and
$$A^- \leq \sum_{i\geq 1} (\sum_{n\geq 1} a_n^2 (-\log_2a_n)_i^2)^\frac12 \leq A^+.$$
$$B^- \leq \sum_{i\geq 1}(\sum_{\underset{2^{-2^{i+1}} <a_n\leq 2^{-2^i}}\to{n\geq 1}} a_n^2 \log_2^2 a_n)^\frac12 \leq B^+,$$
for
$$A^- =\sum_{i\geq 1} 2^i ||u_{i+1} + u_{i+2}+\dots||,$$
$$A^+ =\sum_{i\geq 1} 2^{i} ||u_i +u_{i+1}+\dots||,$$
$$B^- =\sum_{i\geq 1} 2^i ||u_i||,$$
$$B^+ =\sum_{i\geq 1} 2^{i+1} ||u_i||.$$

Then $A^- <\infty \iff A^+ <\infty$, $B^- <\infty \iff B^+ <\infty$, and
$$B^-\leq A^+ \leq \sum_{i\geq 1} 2^i(||u_i||+||u_{i+1}||+\dots) =\sum_{j\geq 1}(2+\dots +2^j)||u_j||\leq B^+.$$

\baselineskip-12pt
\hfill $\square$
\enddemo

\baselineskip18pt

Let us observe, that all conditions $(\alpha), (\beta), (\gamma)$ can be formulated using information function of a partition of some interval.
Once more, assume for simplicity that $\sum |a_n|^2 =1$ and denote $B =\{\sum_{m\geq n} |a_n|^2; n\geq 1\}$.  Let $I_B$: $(1, 1]\to\Bbb R^+$ be given by
$$I_B(t) =-\log_2(\beta -\alpha) \quad\text{for } \ t\in (\alpha, \beta]$$
with $[\alpha, \beta]\cap B =\{\alpha, \beta\}.$

\medskip
Let us fix the following

\medskip
\flushpar{\bf 1.6. Notation.} We write
$$f_0 = f\wedge 2,\quad f_i = (f\wedge 2^{i+1}) - (f\wedge 2^i), \quad i\geq 1$$
 for any positive function $f.$

Gist of the matter given in Theorems 1.3, 1.4 is contained in

\proclaim{1.7. Theorem} A.  For any decreasing sequence $a_n\searrow 0$, $\sum a_n^2 =1$, the condition (*) is equivalent to
$$||I_B|| <\infty.\tag{$\alpha_1$}$$

B. For any $a_n >0$, $\sum a_n^2 =1$, the conditions (**),
$$\sum_{i\geq 1} ||I_B 1_{(2^i\leq H_B <2^{i+1})}|| < \infty,\tag{$\beta_1$}$$
$$\sum_{i\geq 1} ||(I_B)_i||< \infty \tag{$\gamma_1$}$$
are equivalent.

\endproclaim

The equivalence of $(\beta_1)$ and $(\gamma_1)$ is given by 1.5.

For a not  necessarily decreasing,  sequence $a_n\geq 0$,  condition $(\alpha_1)$ is too weak and $(\gamma_1)$ is too strong to characterize the phenomenon (*). It turns out that 
$V I_B <\infty$ is a proper, intermediate, condition. Obviously we use $\log_3$ instead $\log_2$ and $h_B$ instead $I_B$ for technical reasons only.

\bigskip

\flushpar{\bf 2. Fundamental lemma for construction of divergent orthogonal series.} In this section we point out that some orthogonal sequences $\phi_n$, $1\leq n\leq N$, in $L_2(\Bbb R)$ satisfying $\sum_{1\leq n\leq N} ||\phi_n||^2 =1$, have majorants
$$M = \max_{1\leq n\leq N} (\phi_1+\dots +\phi_n)$$
of logarithmic magnitude:
$$\int_\Bbb R M^2 d\lambda \geq k \log^2 N$$
 for some constant $k >0$. We should describe this phenomenon in a specific way suitable for further, rather  complicated calculations. Special properties of the final sum $\phi_1 +\dots +\phi_N$ will also be needed. That is why we don't use the classical constructions based on properties of Hilbert matrices (cf. [2]). Our Lemma 2.1. is obtained by the use of tertiary expansions of numbers and improves some results in [4].

We fix some notation used throughout  the  section.
For $x\in [0, 1)$, we write
$$x = \frac{x_1}{3} +\frac{x_2}{3^2}+\dots,\quad x_k =0,1,2, \ x_k \nrightarrow 2.$$
Each $x_k$ is identified with a function $x_k(x) = x_k$ on $[0, 1)$. 
We use the probabilistic notation where, for example, $(x_k =1) = \{x: x_k(x)= 1\}$.
For a fixed $k$ and any $n =0,\dots, 3^k-1$, we write
$$n = n_1 3^{k-1} + \dots +n_k 3^0,\quad n_k =0,1,2,$$
In the proof of our fundamental lemma we also use the special convention that
$$\aligned &\widehat 0 =0, \quad\widehat 1 =1,\quad \widehat 2=-1,\\
&\widehat{3p+m} =\widehat m \ \text{for } p\in\Bbb Z, \ m =0,1,2.\endaligned\tag0$$
 In other words
$$\widehat m \equiv m \ (\text{mod}\ 3),\quad \widehat m =-1, 0, 1, \quad\text{for } m\in\Bbb Z.$$

\proclaim{2.1. Fundamental Lemma} Let $k =1,2,\dots$ be fixed. For any function $\chi\in L_2(\Bbb R)$, $\chi 1_{[0, 1)}=0$, $||\chi||^2 =1$, there exist functions $\phi_n$, $0\leq n < 3^k$, satisfying
$$\phi_n  \ \text{are mutually orthogonal}, \ ||\phi_n||^2 = \frac{3}{3^k};\tag1$$
$$\langle \phi_n, 1_{[0, 1)}\rangle =0,\quad \sum_{0\leq n < 3^k} \phi_n 1_{[0, 1)} =0;\tag2$$
$$\phi_n 1_{\Bbb R\setminus [0, 1)} = \frac{\sqrt{3}}{3^k} \chi,\tag3$$
in particular $\sum_{0\leq n < 3^k} \phi_n = \sqrt{3} \chi$;
$$\max_{0\leq n <3^k} (\phi_0+\dots +\phi_n) =\sum_{1\leq l\leq k} 1_{(x_l =1)}.\tag4$$
Moreover,
$$\phi_0+\dots +\phi_{n-1} =\sum_{1\leq l\leq k} 1_{(x_l =1)} \tag5$$
and
$$\phi_n=0\tag6$$
for  $x\in [\frac{n}{3^k}, \frac{n+1}{3^k}),$
$0 \leq n < 3^k$.
\endproclaim

\demo{Proof} It is enough to take
$$\phi_n =\frac{1}{3^k} (\sqrt{3} \chi + 3 \widehat{x_1-n_1} + 3^2 1_{(x_1 =n_1)} \widehat{x_2-n_2}+\dots + 3^k 1_{(x_1=n_1,\dots x_{k-1} =n_{k-1})}\widehat{x_k-n_k}).$$
To obtain (1) and to make the last formula more familiar, we start with some properties of functions $1_{(x_1-n_1,\dots, x_{l-1}=n_{l-1})} \widehat{x_l-n_l}$, $1\leq l\leq k$. According to our notation, let  $n' =n'_1 3^{k-1}+\dots + n'_k 3^0$, $n'_l =0,1,2$. For symbol '$\widehat{ \ }$' given by (0), we have
$$||\widehat{x_1 -n_1}||^2 =\frac{2}{3}, \quad \langle\widehat{x_1-n_1}, \widehat{x_1-n'_1}\rangle = -\frac{1}{3} $$
for $n_1\ne n'_1$ and
$$||1_{(x_1 =n_1,\dots, x_{l-1} = n_{l-1})} \widehat{x_l-n_l}||^2 = \frac{2}{3^l},$$
$$\langle 1_{(x_1= n_1,\dots, x_{l-1} =n_{l-1})} \widehat{x_l-n_l}, 1_{(x_1= n_1,\dots, x_{l-1} =n_{l-1})} \widehat{x_l -n'_l}\rangle = -\frac{1}{3^l}$$
for $n_l \ne n'_l$, $1\leq l\leq k$.
Moreover $\langle\widehat{x_1 -n_1}, 1_{[0, 1)}\rangle =0$, 
 $\langle \widehat{x_l-n_l}, 1_{(x_1=n_1,\dots, x_{l-1} =n_{l-1})} \rangle =0$ and
$$\langle 1_{(x_1=n_1,\dots, x_{l-1} =n_{l-1})} \widehat{x_l -n_l}, 1_{(x_1=n'_1,\dots, x_{l'-1} =n'_{l'-1})} \widehat{x_{l'}-n_{l'}}\rangle =0$$
for $l\ne l'$. The orthogonality
$$\langle 1_{(x_1=n_1,\dots, x_{l-1} =n_{l-1})} \widehat{x_l-n_l},  1_{(x_1=n'_1,\dots, x_{l-1} =n'_{l-1})} \widehat{x_l-n'_l}\rangle=0$$
for $(n_1,\dots, n_{l-1})\ne (n'_1,\dots, n'_{l-1})$ is also obvious.

Thus, for $n\ne n'$ with
$(n_1,\dots, n_{l-1}) = (n'_1,\dots, n'_{l-1})$, $n_l\ne n'_l$, we have
$$\multline\langle\phi_n, \phi_{n'}\rangle =\frac{1}{3^{2k}}(||\sqrt{3}\chi||^2 \\
+ 3^2 ||\widehat{x_1-n_1}||^2 +\dots
 +3^{2l-2}||1_{(x_1=n_1,\dots, x_{l-2} =n_{l-2})} \widehat{x_{l-1}-n_{l-1}}||^2\\
+ 3^{2l} \langle 1_{(x_1=n_1,\dots, x_{l-1} 
=n_{l-1})}
 \widehat{x_l-n_l}, 1_{(x_1=n_1,\dots, x_{l-1}=n_{l-1})} \widehat{x_l-n'_l}\rangle)\\
= \frac{1}{3^{2k}} (3 +3^2\cdot \frac{2}{3}+\dots +3^{2l-2}\frac{2}{3^{l-1}} -3^{2l} \frac{1}{3^l}) =0,\endmultline$$
and analogously
$$||\phi_n||^2 =\frac{1}{3^{2k}} (3+3^2\cdot \frac{2}{3}+\dots + 3^{2k} \frac{2}{3^k})=\frac{3}{3^k}.$$
Condition (1) is thus proved.

Conditions (2) and (3) are obvious, by the definition of $\phi_n$ and because
$$\sum_{n_1=0,1,2} \widehat{x_1-n_1} =0,$$
$$\sum_{n_l=0,1,2} 1_{(x_1=n_1,\dots, x_{l-1}=n_{l_1})} \widehat{x_l-n_l} =0$$
for any fixed $l$, $n_1,\dots, n_{l-1}.$

To obtain (4), it is enough to observe that
$$1_{(x_1=m_1)} \widehat{x_1-m_1} \equiv 0 \ \text{for } m_1=0,1,2.$$
Thus
$$1_{(x_1 =n_1,\dots,x_l =n_l)} \widehat{x_l-n_l}\equiv 0$$
and 
$$1_{[\frac{n}{3^k}, \frac{n+1}{3^k})} \widehat{x_l-n_l} =1_{(x_1=n_1,\dots, x_k=n_k)} \widehat{x_l-n_l} \equiv 0.$$

The proof of (4) and (5) is more difficult. We give it with details. We list here the basic properties of the symbol '\ $\widehat{}$ '. We have, for $\lambda\in [0, 1],$
$$\aligned \lambda \widehat{x_1-0} \leq 1_{(x_1=1)} \quad &\text{for } x\in [0, 1),\\
\lambda \widehat{x_1-0} =1_{(x_1=1)} \quad &\text{for } x\in [0, \frac{1}{3});\\
{}&{}\\
\widehat{x_1-0} +\lambda\widehat{x_1-1} \leq 1_{(x_1=1)} \quad &\text{for } x\in [0, 1),\\
\widehat{x_1-0} +\lambda\widehat{x_1-1} = 1_{(x_1=1)} \quad &\text{for } x\in [\frac{1}{3}, \frac{2}{3});\\
{}&{}\\
\widehat{x_1-0} +\widehat{x_1-1} +\lambda\widehat{x_1-2} \leq 1_{(x_1=1)} \quad &\text{for } x\in [0, 1),\\
\widehat{x_1-0} +\widehat{x_1-1} +\lambda\widehat{x_1-2} = 1_{(x_1=1)} \quad &\text{for } x\in [\frac{2}{3}, 1);\\
{}&{}\\
\widehat{x_1-0} + \widehat{x_1-1} + \widehat{x_1 -2} =0 \quad &\text{for } x\in [0, 1).\endaligned$$
More generally, for $0\leq l\leq k$, $\lambda\in [0, 1]$,
$$\aligned\lambda\widehat{x_l-0} \leq 1_{(x_l=1)} \quad &\text{for } x\in [0, 1),\\
\lambda \widehat{x_l-0} =1_{(x_l=1)} \quad &\text{for } x\in (x_l =0);\\
{}&{}\\
\widehat{x_l-0} +\lambda\widehat{x_l-1} \leq 1_{(x_l=1)} \quad &\text{for } x\in [0, 1),\\
\widehat{x_l-0} +\lambda\widehat{x_l-1} =1_{(x_l=1)}\quad &\text{for } x\in (x_l=1);\\
{}&{}\\
\widehat{x_l-0} +\widehat{x_l-1} +\lambda \widehat{x_l-2} \leq 1_{(x_l=1)} \quad &\text{for } x\in [0, 1),\\
\widehat{x_l-0} +\widehat {x_l-1} +\lambda \widehat{x_l-2} =1_{(x_l=1)}\quad &\text{for } x\in (x_l=2);\\
{}&{}\\
\widehat{x_l-0} +\widehat{x_l-1} +\widehat{x_l-2} =0\quad &\text{for } x\in [0, 1).\endaligned$$
For any $m=0,\dots, 3^k-1$, $m =m_1 3^{k-1} + \dots +m_k 3^0$, we write also $m = [m_1,\dots, m_k]$ with $m_l =0,1,2$, in particular, $[2,\dots, 2]= 3^k-1$. Then, for $1\leq l\leq k$, $x\in [0, 1)$,
$$\multline\sum_{[n_1,\dots, n_{l-1}, 0,\dots, 0]\leq m\leq [n_1,\dots, n_{l-1}, 0, n_{l+1},\dots, n_k]} \widehat{x_l-m_l}\\
= (n_{l+1} 3^{k-l-1} +\dots + n_k 3^0) \widehat{x_l-0} \leq 3^{k-l} 1_{(x_l=1)}\\
\text{with equality for } \ x\in (x_l=0), \endmultline\tag7$$
 and similarly
$$\multline\sum_{[n_1,\dots, n_{l-1}, 0,\dots, 0]\leq m\leq [n_1,\dots, n_{l-1}, 1, n_{l+1}, \dots, n_k]} \widehat{x_l-m_l}
\leq 3^{k-l} 1_{(x_l=1)}\\
 \text{with equality for } x\in (x_l=1),\endmultline\tag8$$
$$\multline\sum_{[n_1,\dots, n_{l-1}, 0,\dots, 0]\leq m\leq [n_1,\dots, n_{l-1}, 2, n_{l+1},\dots, n_k]} \widehat{x_l-m_l} \leq 3^{k-l} 1_{(x_l=1)} \\
\text{with equality for } x\in (x_l =2).\endmultline\tag9$$
 Moreover,
$$\sum_{[n_1,\dots,n_{l-1}, 0,\dots,0]\leq m\leq [n_1,\dots, n_{l-1}, 2,\dots, 2]} \widehat{x_l-m_l} = 3^{k-l}(\widehat{x_l-0} +\widehat{x_l-1} +\widehat{x_l-2}) =0.\tag10$$
Both conditions (4) and (5) are consequences of
$$\multline \sum_{0\leq m <n} 1_{(x_1=m_1,\dots, x_{l-1}=m_{l-1})} \widehat{x_l-m_l}\\
\leq 3^{k-l} 1_{(x_l=1)} \quad\text{with equality for } x\in (x_1=n_1,\dots, x_l =n_l).\endmultline\tag11$$
Indeed, recall that $x\in [\frac{n}{3^k}, \frac{n+1}{3^k})\iff (x_1 =n_1,\dots, x_k=n_k)$.

We prove relations (11)  separately on each  fixed interval $x\in (x_1=n'_1,\dots, x_{l-1} =n'_{l-1})$; thus we prove that
$$\multline \sum_{\underset{[n'_1,\dots, n'_{l-1}, 0,\dots, 0]\leq m\leq [n'_1,\dots, n'_{l-1}, 2,\dots, 2]}\to{0\leq m <n}}
1_{(x_1=m_1,\dots x_{l-1}=m_{l-1})} \widehat{x_l-m_l} \\
\leq 3^{k-l} 1_{(x_l=1)}\quad
\text{with equality if only } x_1=n_1 =n'_1,\dots, x_{l-1} =n_{l-1} = n'_{l-1},  x_l =n_l.\endmultline\tag12$$
We discuss three cases:

1$^o$ Assume that $[n_1,\dots, n_{l-1}, 0,\dots, 0]< [n'_1,\dots, n'_{l-1}, 0,\dots, 0]$. Then we have zero summands in (12) and everything is obvious.

2$^o$ Assume that $[n_1,\dots, n_{l-1}, 0,\dots, 0]> [n'_1,\dots, n'_{l-1}, 0,\dots, 0]$. Then  the sum  can be be written as
$$\sum_{[n'_1,\dots, n'_{l-1}, 0,\dots, 0]\leq m\leq [n'_1,\dots, n'_{l-1}, 2,\dots, 2]} 1_{(x_1=n'_1,\dots, n'_{l-1}=n'_{l-1})} \widehat{x_l-m_l}$$
and it equals $0$ by (10).

3$^o$ Assume that $n_1=n'_1,\dots, n_{l-1} =n'_{l-1}$. Then the sum in (12) equals
 $$\sum_{[n_1,\dots, n_{l-1}, 0,\dots, 0]\leq m < [n_1,\dots, n_{l-1}, n_l\dots, n_k]} 
1_{(x_1=m_1,\dots, m_{l-1}=n'_{l-1})} \widehat{x_l-m_l}.$$
The required relations can be obtained from (7) if $n_l =0$, from (8) if $n_l =1$, and from (9) if $n_l =2$. The proof is finished.
\hfill $\square$
\enddemo

For large $l$ the function $\sum_{1\leq l\leq k} 1_{(x_l=1)}$ approximates $\frac{1}{3}k$. By Bernstein inequality we have in particular

\proclaim{2.2. Lemma} For $x= \frac{x_1}{3}+ \frac{x_2}{3^2}+\dots$, $x_l =0,1,2$, we have
$$\lambda(\sum_{1\leq l\leq k} 1_{(x_l =1)} < \frac{1}{6} k) \leq e^{-k/144}.$$
\endproclaim

\demo{Proof} The left-hand side equals $P(\frac{1}{k} S_k -p < -\epsilon)$ for the Bernoulli random variable $S_k$ with probability of success  $p =\frac{1}{3}$, and for $\epsilon =\frac{1}{6}$. Thus classical inequality $P(\frac{1}{k} S_k -p < -\epsilon) \leq e^{-k \epsilon^2/4}$ (see [1]) can be used.
\enddemo

\bigskip
\flushpar{\bf 3. Consequences of the fundamental lemma for triadic sets}. In the rest of the paper  (excluding Section 6) $B$ is always a set satisfying
$$\aligned &\{0, 1\}\subset B\subset [0, 1],\\
& \sharp B \cap [\alpha, 1] <\infty\quad\text{for any } \alpha >0.\endaligned\tag13$$

\definition{3.1. Definition} We say that a set $B$ satisfying (13) is {\it triadic} if
$$B = \bigcup_{(i,n)\in I} \{n 3^{-2^i}, (n+1) 3^{-2^i}\}$$
for a set $I$ of pairs $(i, n)$ $i\geq 0$, $0\leq n <3^{2^i}$ and
$$(n 3^{-2^i}, (n+1) 3^{-2^i})\cap B\ne \emptyset \quad\text{implies}\quad \{n 3^{-2^i}, (n+1) 3^{-2^i}\}\subset B,$$
for any $i\geq 0$, $0\leq n < 3^{2^i}$. We also assume that $\{0, \frac{1}{3}, \frac{2}{3}, 1\}\subset B$.
\enddefinition

\definition{3.2. Definition} We say that $X: B\to L_2(\Bbb R)$ is an {\it orthogonal process} if $X(0) =0$ and
$$||X(t) - X(s)||^2 = t-s \quad\text{for any } 0\leq s <t \leq 1.$$
\enddefinition

Sometimes we have $X(B) \subset L_2[a, b)$; then $L_2[a, b)$ is identified with the space of functions vanishing on outside of $[a, b)$.

\medskip

\flushpar{\bf 3.3. Basic lemma for a finite triadic set}. There exists a constant $c >0$ such that, for any finite triadic set $B$ and any $y >1$, if $V h_B > c y$ then
$$\lambda([0, 1)\setminus (\max_{t\in B} X(t) >y)) <\frac12$$
for some orthogonal process $X$.

Thus in this section we discuss a {\it fixed finite triadic set} $B$. We need a number of auxiliary lemmas and notations connected with $B$.

\medskip
\definition{ 3.4. Definition} We say that $D$ is simple, and we write $D\in\Cal S^B$, if $D$ is a finite union of closed intervals with end points in $B.$

\enddefinition

\definition{3.5. Definition} Let $D \in \Cal S^B$. We say that $X$ is a {\it simple process} on $D$, and we write $X\in \Cal P_D^B$, if $X: D\cap B\to L_2(\Bbb R)$ and
$$||X(t) - X(s)||^2 = 3\cdot 24^2 \lambda ([s, t]\cap D), \quad s,t,\in D\cap B;$$
$$X(\alpha) =0\quad\text{for }\quad \alpha =\min D.$$
\enddefinition
For any orthogonal process $X$ on $B$ and $D = [\alpha, \beta]$, $\{\alpha, \beta\}\in B$, the renormalized process  $Y(t) = 24\sqrt{3}(X(t) - X(\alpha))$ defined on $D\cap B$ is a simple process on $D.$

\definition{3.6. Definition} For a given $y >0$, $0 <\epsilon <1$, we say that a simple set $D$ is a set of $(\epsilon, y)$ -- complexity for $B$, and we write $D\in \Cal S^B(\epsilon, y)$, if for any interval $[a, b)$ and any $\chi\in L_2(Z)$, $||\chi||=1$, $Z\cap [a,b) =\emptyset$, there exists simple process on $D$ satisfying
$$X(\max D) = 24\sqrt{3 \lambda(D)} \chi,$$
$$\lambda([a, b) \setminus (\max_{t\in D\cap B}  X(t) \geq \frac{y}{\sqrt{b-a}})) <\epsilon(b -a),$$
$$X(t) \in L_2([a, b])\cup Z)\quad\text {for } t\in D\cap B.$$
\enddefinition

\flushpar{\bf 3.7. Example.} By Lemmas 2.1, 2.2 we have
$$[0, 1]\in\Cal S^B(e^{-k/144}, 4k)\tag14$$
 if only $\{m 3^{-k}; 0\leq m\leq 3^k\}\subset B$. Namely, for $[a, b) = [0, 1)$, $\chi\in L_2(Z)$, $||\chi|| =1$, $Z\cap [0, 1) =\emptyset$, it is enough to take
$$\tilde X(m 3^{-k}) =24(\phi_0 +\dots + \phi_{m-1}), \quad 0\leq m\leq 3^k,$$
and then extend $\tilde X$ to a simple  process
$$X: B\to L_2([0, 1)\cup Z).$$

For an arbitrary interval $[a, b)$, it is enough to take $\phi_n(\frac{t-a}{b-a}) /\sqrt{b-a}$ instead of $\phi_n(t)$ (after suitable rearangement of $X$ and $Z$).

It is also obvious that
\medskip
\flushpar{\bf 3.8. Example.} $[\alpha, \beta]\in \Cal S^B(e^{-k/144}, 4k\sqrt{\beta-\alpha})$
if only $\{\alpha +m(\beta -\alpha) 3^{-k}$;  $0\leq m\leq 3^k\}\subset B$.

Moreover, let $D=\bigcup_{0\leq n <3^{k}}\delta_n$ with $\delta_n$ being closed intervals
from  $\Cal S^B$ (with end points in $B$), with mutually disjoint interiors and  with $\lambda(\delta_n) =\eta$. Then
$$D\in\Cal S^B(e^{-k/144},\quad 4k\sqrt{\lambda(D)}).\tag15$$

\flushpar{\bf 3.9. Remark.} The multiplier 4 in formulas (14), (15) and their just comming generalizations is suitable in further, more complicated considerations. This is the only reason of using the  strangely looking  constant $3\cdot 24^2$ in Definition 3.5.
\medskip
The goal of this section is to show that, for some $C >0$
$$||V h_B|| > C y \quad\text{implies} \quad [0, 1]\in \Cal S^B(\frac12, y)\tag16$$
for any $y\geq 1$. Then  Lemma 3.3 is proved.
\medskip

We need more delicate consequences of Lemmas 2.1 and 2.2. Let us observe that
$$\Cal S^B(\epsilon, y)\subset \Cal S^B(\epsilon_1, y_1)$$
for $\epsilon_1 \geq \epsilon$, $y_1\leq y$, and

\proclaim{3.10. Lemma} If sets $D_1,\dots, D_L$ have mutually disjoint interiors and $D_l\in\Cal S^B(\epsilon, y_l)$, $1\leq l\leq L$, then
$$\bigcup_{1\leq l\leq L} D_l\in\Cal S^B(\epsilon, \sqrt{\sum_{1\leq l\leq L} y_l^2}).$$
\endproclaim

\demo{Proof} For any interval $[a, b)$ and any function $\chi\in L_2(Z)$, $||\chi|| =1$, let us take a partition
$$[a, b) = [a_1, b_1)\cup\dots \cup [a_L, b_L)$$
with disjoint intervals satisfying
$$\frac{b_l-a_l}{b-a} = y_l^2/\sum_{1\leq l'\leq L} y^2_{l'}\quad \text{for } 1\leq l\leq L,$$
and a decomposition
$$\chi =\sum_{1\leq l\leq L} \frac{\sqrt{\lambda(D_l)}}{\sqrt{\lambda(\bigcup_{l'} D_{l'})}} \chi_l$$
with some orthonormal system $(\chi_1,\dots, \chi_L)$ in $L_2(Z)$. By assumption, there exist simple processes $X_l\in \Cal P^B_{D_l}$ such that
$$\lambda([a_l, b_l)\setminus (\max_{t\in D\cap B} X(t) \geq \frac{y_l}{\sqrt{b_l-a_l}}))<\epsilon (b_l-a_l),$$
$$X_l(\max D_l) = 24\sqrt{3} \sqrt{\lambda(D_l)} \chi_l,$$
$$X_l(t)\in L_2([a_l, b_l)\cup Z),$$
for $1\leq l\leq L$. It is enough to take
$$X(t) =\sum_{1\leq l\leq L} X_l(\max(D_l\cap [0, t]))$$
for $t\in D\cap B$.
\hfill $\square$
\enddemo

\proclaim{3.11. Lemma} Let $D =\bigcup_{0\leq n <3^k} \delta_n$, $\delta_n$ being closed intervals from $\Cal S^B$, with mutually disjoint interiors. If $\lambda(\delta_n) =\eta$ and $\delta_n\in \Cal S^B(\epsilon, y)$, $0\leq n <3^k$, then
$$D\subset \Cal S^B(\epsilon +e^{-k/144}, 3^{k/2}y+ 4k\sqrt{\lambda(D)}).$$

\endproclaim

\demo{Proof} One can assume, by a suitable change of notations, that $\delta_n =[\alpha_n, \beta_n]$, $\beta_n\leq \alpha_{n+1}$, $0\leq n <3^k-1$. Take for simplicity $[a, b) =[0, 1)$ and $\chi\in L_2(Z)$, $||\chi|| =1$, $[0, 1)\cap Z =\emptyset$. By Lemma 2.1, there exist functions $\phi_n$ satisfying
$$||\phi_n||^2 =3\cdot 24^2\eta,\quad \phi_n\in L_2[0, 1)\cup Z,$$
and also
$$\phi_n \bot L_2 [n 3^{-k}, (n+1) 3^{-k}),$$
for  $0\leq n < 3^k$. Moreover
$$\sum \phi_n = 24 \sqrt{3\cdot 3^k\eta} \chi = 24\sqrt{3\lambda(D)} \chi,$$
$$(\phi_0+\dots +\phi_{n-1})(x) =\sum_{1\leq i\leq k} 1_{(x_i=1)}(x) \quad\text{for } x\in [n 3^{-k}, (n+1) 3^{-k}),$$
for any $0\leq n <3^k$.

For a given $0\leq n <3^k$, we use our assumption $\delta_n\in \Cal S^B(\epsilon, y_n)$ with $[a, b)= [n 3^{-k}, (n+1) 3^{-k})$, $\chi =\phi_n/||\phi_n||\in L_2([0, n 3^{-k})\cup [(n+1) 3^{-k}, 1)\cup Z)$ and we obtain a simple process $X_n$ on $\delta_n\cap B$, $X_n\in\Cal P^B_{\delta_n}$, $X_n(\alpha_n) =0$, $X_n(\beta_n) =\phi_n$,
$$\lambda([ n 3^{-k}, (n+1) 3^{-k})\setminus (\max_{t\in \delta_n\cap B} X_n(t) \geq 4\sqrt{3^k} y)) <\epsilon 3^{-k},$$
$$X_n(t)\in L_2, ([0, 1)\cup Z).$$
Obviously, it is enough to take $X$ satisfying
$$X(t) =\phi_0+\dots +\phi_{n-1} +X_n(t)$$
for $t\in (\alpha_n, \beta_n]$, $0\leq n <3^k.$
\hfill $\square$

\enddemo

Let us reformulate Lemmas 3.10, 3.11 in  a more useful way

\proclaim{3.12. Corollary} A. If $D_1,\dots, D_L$ have mutually disjoint interiors, $D = \bigcup_l D_l$ and
$$D_l \in\Cal S^B(\epsilon, ||h 1_{D_l}||),\quad 1\leq l\leq L,$$
for some function $h\in L_2(0, 1]$, then
$$D\in \Cal S^B(\epsilon, ||h 1_D||).$$

B. If $\delta_1,\dots , \delta_{3^k}$ are closed intervals of the same length with mutually disjoint interiors, $D= \bigcup_l\delta_l$ and
$$\delta_l\in\Cal S^B(\epsilon, ||a 1_{\delta_l}||),\quad 1\leq l\leq 3^k,$$
$a\geq 0$, then
$$D\in\Cal S^B(\epsilon + e^{-k/144}, ||(a + 4k) 1_D||).$$
\endproclaim

Now we present some properties of information  function $h_B$ defined in Section 1.
Let us put
$$\delta_n^j = (n 3^{-2^j}, (n+1) 3^{-2^j}],\quad 0\leq n <3^{2^j}, \quad j\geq 0.\tag17$$

\definition{3.13. Definition} We say that a bounded Borel function $h$ on $(0, 1]$ is {\it triadic} if $h\geq 1$ and
$$\delta_n^j \cap (h\geq 2^j)\ne \emptyset$$
implies
$$\delta_n^j \subset (h\geq 2^j),$$
for any $j\geq 0$, $0\leq n <3^{2^j}.$

\enddefinition

Thus obviously

\proclaim{3.14. Lemma } For any triadic set $B$, the information function $h_B$ is triadic.

\endproclaim

\definition{3.15. Definition} A triadic function $h$ is of type $j$ $(h\in \Cal T_j$  in symbols), if

$1^o$ $h$ is constant on each $\delta_n^{j+1}$, $0\leq n <3^{2^{j+1}}$;

$2^o$ for each value $h(t) <2^{j+1}$ we have $h(t) = 2^{j'}$ for some $0\leq j'\leq j.$

In particular, for any function $h\in \Cal T_j$ and any $\delta_m^j\subset (h\geq 2^j)$, $0\leq m <3^{2^j}$, a representation
$$(h-2^j) 1_{\delta_m^j} =\sum_{1\leq k\leq K} (2^j+ a_k) 1_{\delta^{j+1}_{n(k)}}\quad
\text{with $a_k\geq 0$}\tag18$$
 is possible, and
$\delta_m^j\cap (h\geq 2^{j+1}) =\bigcup_{1\leq k\leq K} \delta^{j+1}_{n(k)}.$

\enddefinition

Denoting
$$\underline{a} = 2^j\quad\text{for } \ 2^j \leq a <2^{j+1}, \ j\geq 0,\tag19$$
for any $a\geq 1$, we have

\proclaim{3.16. Lemma } For any triadic set $B$ with information function $h_B\leq i+1$, 
$\underline{h_B}$ is of type $i$, $\underline{h_B} \in\Cal T_i.$
\endproclaim

\definition{3.17. Definition} For any positive sequences $(a_k), (b_k)$, $1\leq k\leq K$, we write
$$(b_k) \prec_j (a_k),\quad j\geq 1,$$ 
if, for some mutually disjoint classes of indices $I_s\subset \{1,\dots, K\}$, $\sharp I_s =3^{2^{j-1}}$ for $s\in S$, we have
$$b_k \leq 2^{j+1} +\min_{l\in I_s} a_l\quad\text{for } k\in I_s, \ s\in S,$$
$$b_k=a_k \quad\text{for} \quad k\notin \bigcup_{s\in S} I_s.$$
\enddefinition

\definition{3.18 Definition} We say that a (nonlinear) operator $U$ is of type $j$ if  $U$
is defined on $\Cal T_j,$ and

$1^o$ $U h\wedge 2^j = h\wedge 2^j;$

$2^o$ for $\delta_m^j \subset (h\geq 2^j)$, we have $\delta_m^j\subset (U h\geq 2^j)$ and
$$(U h -2^j) 1_{\delta_m^j} =\sum_{1\leq k\leq K} b_k 1_{\delta_{n(k)}^{j+1}}$$
for some $(b_k) \prec_j(a_k)$,
according to representation (18).

\enddefinition

\proclaim{3.19. Lemma } If we have $h\in \Cal T_j$, then $V_j h\in\Cal  T_{j-1}$ and $V_jU_jh\in\Cal T_{j-1}$ for any operation $U_j$ of type $j.$
\endproclaim

The notion of type $j$ operators is useful, because it is natural to describe the complexity of $[0, 1]$ for a triadic set $B$ by a norm $||U_8V_8\dots U_iV_i h_B-2^8||$ with $U_j$ of type $j$ (at first). It is done in the following two lemmas. The convention $\Cal S^B(e, y) =\Cal S^B(\epsilon, 0) = \ \text{Borel} [0, 1]$ for $e\geq 1$ is natural.

\proclaim{3.20. Lemma} Let
$$h\wedge 2^{j+1} = \underline{h_B} \wedge 2^{j+1}$$
for some $h\in \Cal T_j$ and some triadic set $B$. If we have
$$\delta_n^{j+1} \in \Cal S^B(\epsilon_j, ||(h-2^{j+1})^+ 1_{\delta_n^{j+1}}||)$$
for any $0 \leq n < 3^{2^{j+1}}$, then
$$\delta_m^j\in\Cal S^B(\epsilon_j +e^{-2^{j-1}/144}, ||(V_jU_j h-2^j)^+ 1_{\delta_m^j})||)$$
for $0\leq m <3^{2^j}$ and for any operator $U_j$ of type $j$. Moreover,
$$V_jU_j h\in \Cal T_{j-1}, \quad V_jU_j h\wedge 2^j = \underline{h_B} \wedge 2^j.\tag20$$
\endproclaim

\demo{Proof} Fix $\delta_m^j\in (\underline{h_B} \geq 2^j)$. For $I_s$, $s\in S_m$, defined as in 3.17, 3.18, we have
$$\bigcup_{n\in I_s} \delta_n^{j+1}\in\Cal S(\epsilon_j +e^{-2^{j-1}/144}, ||(m_s +4\cdot 2^{i-1}) 1_{\bigcup_{n\in I_s} \delta^{j+1}_n}||)$$
with $m_s =\min \{h(t)-2^{j+1}; t\in \bigcup_{n\in I_s}\delta_n^{i+1}\}$, by 3.12.B.

Thus
$$\bigcup_{n\in I_s} \delta_n^{j+1}\in\Cal S(\epsilon_j +e^{-2^{j-1}/144}, ||(U h-2^j)^+ 1_{\bigcup_{n\in I_s^j}\delta_n^{j+1}}||).$$
Then we use 3.12.A, with
$$ \{D_1,\dots, D_L\} = \{\bigcup_{n\in I_s} \delta_n^{j+1}; s\in S_m\}
\cup \{\delta_n^{j+1}; n\notin \bigcup_{s\in S_m} I_s, \delta_n^{j+1}\subset \delta_m^j\},$$
to obtain
$$\delta_m^j\in \Cal S^B(\epsilon_j+ e^{-2^{j-1}/144}, ||(U_j h-2^j)^+ 1_{\delta_m^j}||).$$
The equality
$$||(U_j h-2^j)^+ 1_{\delta_m^j}|| = ||(V_jU_j h-2^j)^+ 1_{\delta_m^j}||$$
and relations (20) are obvious.
\hfill $\square$ 
\enddemo

\proclaim{3.21. Lemma} Let $U_j$ be an operator of type $j$ for $8\leq j\leq i$, and let $B$ be a triadic set with $h_B \leq 2^{i+1}.$
Then
$$[0, 1]\in\Cal S^B(\frac12, ||(V_8U_8\dots V_iU_i \underline{h_B}-2^8)^+||).$$
\endproclaim
\demo{Proof} The assumptions of the previous Lemma, for $j=i$, $h=\underline{h_B}$ and $\epsilon_i=0$, are obviously satisfied and
$$\delta_m^i\in\Cal S^B(e^{-2^{i-1}/144}, ||(V_i U_i \underline{h_B} -2^i) 1_{\delta_m^i}||)$$
for any $\delta_m^i\subset (h_B \geq 2^i)$.
Then, by backward induction,
$$\delta_m^j\in\Cal S^B(e^{-2^{j-1}/144} +\dots + e^{-2^{i-1}/144}, ||(V_jU_j\dots V_iU_i \underline{h_B}-2^j)^+ 1_{\delta_m^j}||)$$
for any $\delta_m^j\subset (h_B\geq 2^j).$

It is enough to observe that $e^{-2^{8-1}/144}+\dots + e^{-2^{i-1}/144} <\frac12$, and use once more 3.19.A with
$$\{D_1,\dots, D_L\} =\{\delta_m^8; 0\leq m <2^8, \delta_m^8\subset (\underline{h_B} \geq 2^8)\}.$$
\vskip-12pt
\hfill $\square$
\enddemo

\vskip18pt

Let us stress that Lemma 3.21 containes the main idea of the proof of implication (16) (and of Lemma 3.3). We only need the implication
$$\multline V h_B \geq C y,\quad y >1, \quad\text{implies } \ 
||V_8U_8\dots V_iU_i \underline{h_B} -2^8)^+||\geq y\\
\text{for suitably chosen $i\geq 1$ and $U_j$ of type $j$, }  8\leq j\leq i.\endmultline\tag21$$

Implication (21) is a consequence of auxliary  Lemmas 3.24, 3.25, 3.26, 3.27.
Namyly, for same $U_j$ of type $j$, the operation $V_j U_j$ is $j$-triadic in the following sense.

\definition{3.22. Definition}  We say that $W_j: \Cal T_j\to\Cal T_{j-1}$, $j \geq 1$, is a $j$-{\it triadic} operation if, for any $h\in\Cal T_j$,
$$W_j h= V_j h- p -q$$
for some positive functions $p, q$ with
$$p+ q \leq (V_j h -2^j)^+,$$
$$p = ||p||_j \leq 2^{-j},\quad q= ||q||_j \leq 2^{-j}(V_j h - 2^j)^+.$$
\enddefinition

Thus in particular $W_j h\wedge 2^j =h\wedge 2^j$ and
$$p+ q,  2^jq\leq ||(h -2^j)^+||_j$$
(cf. Lemma 3.19).

We start with some elementary properties of relation $\prec_i$ defined in 3.14.

\proclaim{3.23. Lemma} For any positive sequence $(a_k)$, $1\leq k\leq K$, and for $i\geq 1$, there exists a sequence $(b_k)$ satisfying $(b_k)\prec_i (a_k)$ and
$$b_k = a_k \quad\text{or } \ a_k +2^i \ \text{for any } 1\leq k\leq K,$$
$$a_k +2^i - b_k = c_k+d_k$$
with some $c_k, d_k \geq 0$,
$$\sum_{1\leq k\leq K} c_k^2 \leq 3^{2^{i-1}}\cdot 2^{2i}\cdot(2^i+1),$$
$$d_k \leq 2^{-i} (a_k+2^i),\quad 1\leq k\leq K.$$

\endproclaim

\demo{Proof} Changing notation, if necessary, one can assume that $(a_k)$, $1\leq k\leq K$, is an increasing sequence. Denote $L =\max\{k; a_k\leq 2^{2i}\}$, $\nu =3^{2^{i-1}}$ and
$$I_1 = \{1,\dots \nu\},\dots,  I_t  =\{(t-1)\nu+1,\dots t\nu\}$$
with $t$ defined by $t\nu \leq L < (t+1)\nu$. Then $\{I_s; s\in S\}$ can be defined by taking the set of indices
$$S =\{s = 1,\dots, t; m_s\geq M_s - 2^i\}$$
with $m_s =\min_{k\in I_s} a_k$, $M_s =\max_{k\in I_s} a_k$.

Take
$$b_k =a_k+2^i\quad \text{for } \ k\in\bigcup_{s\in S} I_s,$$
$$b_k =a_k\quad\text{otherwise}.$$
Then  inequalities
$$b_k =a_k+2^i\leq M_s +2^i \leq m_s+ 2^{i+1} $$
are valid for $k\in I_s$, $s\in S$ and $(b_k)\prec_j(a_k)$ (cf. 3.17). Obviously,
$$2^i\leq 2^{-i}(a_k +2^i)$$ 
for $a_k\geq 2^{2i}$, in particular for $k > (t +1)\nu$.

Note that,  by definition of $L$, $t+1 -\sharp S \leq  \frac{2^{2i}}{2^i} +1 =2^i +1$, and
$$ \sum_{\underset{k\notin\bigcup_{s\in S} I_s}\to{1\leq k\leq (t+1)\nu}}
(2^i)^2 = ((t+1)\nu - \sharp \bigcup_{s\in S} I_s) \cdot 2^{2i} \leq 3^{2^{i-1}}2^{2i}(2^i +1).$$
It is enough to take
$$\aligned c_k & = 2^i\quad\text{for }\quad 1\leq k \leq (t+1)\nu, \ k\notin \bigcup_{s\in S} I_s,\\
c_k &=0 \  \ \text{otherwise};\endaligned$$
$$\aligned d_k & = 2^i\quad\text{for } \quad k>(t+1)\nu\\
d_k &= 0 \ \ \text{otherwise}.\endaligned$$

\vskip-12pt
\hfill $\square$
\enddemo

\vskip18pt

Now we can pass the first and  main step of the proof of implication (21).

\proclaim{3.24. Lemma} For any $j\geq 5$, there exists an operator $U$ of type $j$ such that $W_j = V_j U$ is $j$-triadic, according to 3.22.
\endproclaim

\demo{Proof} Fix $h_j\in \Cal T_j$. For any $\delta_m^j\subset (h_j\geq 2^j)$ we have a representation (cf. (18)),
$$(h_j -2^j)^+ 1_{\delta_m^j} =\sum_{1\leq k\leq K_m} (a_k+2^j)\delta^{j+1}_{n(m, k)},  \  a_k\geq 0.$$
Let
$$(Uh_j -2^j)^+ 1_{\delta^j_m} =\sum_{1\leq k\leq K_m} b_k \delta^{j+1}_{n(m,k)}$$
with $b_k = a_k+2^j  - c_k^m -d_k^m,$
$$\sum_{1\leq k\leq K_m} (c_k^m)^2 \leq 3^{2^{j-1}}\cdot 2^{2j}\cdot(2^j+1),$$
$$d_k^m\leq 2^{-j} (a_k+2^j).$$
Thus $(U h-2^j)^+ = (h-2^j)^+ - f_j -g_j$ (and also $U h = h-f_j-g_j))$ for  
$$f_j =
\sum_{0\leq m\leq 3^{2^j}} \sum_{1\leq  k\leq K_m} c_k^m 1_{\delta_{n(m, k)}^{j+1}}$$
$$g_j =\sum_{0\leq m\leq 3^{2^j}} \sum_{1\leq k\leq K_m} d_k^m 1_{\delta^{j+1}_{n(m,k)}}.$$
Obviously,
$$g_j \leq 2^{-j}(h_j -2^{j})^+,$$
$$||f_j||_j^2\leq  [3^{2^{j-1}} 2^{2j}(2^j +1)]\cdot 3^{-2^{j+1}}/3^{-2^j}\leq 2^{-j}\quad\text{on } \ \delta_m^j$$
if only $j\geq 5.$

Now, $V_jU_j h_j = h_j\wedge 2^j +||(h_j -2^j)^+ - f_j-g_j||_j\geq h_j\wedge 2^j + ||(h_j -2^j)^+||_j - ||f_j||_j - ||g_j||_j$ and we can choose the required functions $p\leq ||f_j||_j$, $q_j\leq ||g_j||_j$.
\hfill $\square$
\enddemo

\proclaim{3.25. Lemma} For any $j$-triadic operations $W_j$, $8\leq j\leq i$, and any $h\in\Cal T_i$, we have
$$||W_8\dots W_ih -2^8||\geq (1-2^{-7})||V_8\dots V_ih -2^8||- 2^{-7}.$$
\endproclaim

\demo{Proof} We have
$$W_j\dots W_ih -V_j W_{j+1}\dots W_ih = p_j +q_j,$$
$$p_j = ||p_j||_j \leq 2^{-j},$$
$$\multline q_j = ||q_j||_j\leq 2^{-j}||(W_{j+1}\dots W_i h-2^j)^+||_j\\
\leq 2^{-j} ||(V_{j+1}\dots V_i h-2^j)^+||_j = 2^{-j}(V_j\dots V_i h-2^j)^+\endmultline$$
and, by a backward induction on $k$,
$$||V_k\dots V_jh_j - V_k\dots V_j(h_j-f)||_k \leq ||f||_k$$
for any triadic function $h_j$, $h_j\wedge 2^j \equiv (h_j -f)\wedge 2^j.$

Thus
$$\aligned ||V_k\dots V_{j-1} W_j\dots & W_ih - V_k\dots V_j W_{j+1}\dots W_i h||_k\\
& \leq ||p_j||_k + ||g_j||_k\\
& \leq 2^{-j} + 2^{-j} ||(V_j\dots V_i h - 2^j)^+||_k\\
& \leq 2^{-j} + 2^{-j} (V_k\dots V_i h-2^k),\endaligned$$
in particular
$$\multline ||V_8\dots V_ih -W_8\dots W_ih||_8\\
\leq (2^{-8} +\dots + 2^{-i}) + (2^{-8} +\dots + 2^{-i})(V_8\dots V_i h-2^8).\endmultline$$
\vskip-12pt
\hfill $\square$
\enddemo 

\baselineskip18pt

\proclaim{3.26. Lemma } For any triadic function $h$ on $(0, 1]$, $h\leq 2^{i+1}$, we have
$$||(V_8\dots V_i h- 2^7)^+|| \leq 3||(V_8\dots V_i\underline{h}-2^7)^+||.$$
\endproclaim

\demo{Proof} For any interval $\delta_m^i = (m 3^{-2^i}, (m+1) 3^{-2^i}]\subset (h\geq 2^i)$ we have $V_i\underline{h} \geq 2^i$, $V_i h\leq 2^{i+1}$ on $\delta_m^i$, and
$$||(V_i h- 2^{i-1}) 1_{\delta_m^i}||\leq 3||(V_i \underline{h} -2^{i-1}) 1_{\delta_m^i}||.$$

Assume that $\delta_m^j = (m 3^{-2^j}, (m+1) 3^{-2^j}] \subset (h\geq 2^j)$, and
$$||(V_{j+1}\dots V_i h -2^j) 1_{\delta_n^{j+1}} ||\leq 3 ||(V_{j+1}\dots V_i \underline{h} -2^j) 1_{\delta_n^{j+1}}||$$
for any $\delta_n^{j+1} = (n 3^{-2^{j+1}}, (n+1) 3^{-2^{j+1}}]\subset \delta_m^j\cap (h \geq 2^{j+1})$, and let $I$ be the set of such indices $n$. Then
$$\multline ||(V_j V_{j+1}\dots V_i h-2^{j-1}) 1_{\delta_m^j}|| = ||(2^j-2^{j-1}) 1_{\delta_m^j} || +||(V_{j+1}\dots V_i h -2^j) 1_{\delta_m^j}||\\
\leq ||2^j-2^{j-1} 1_{\delta_m^j}|| + ||(2^{j+1}-2^j) 1_{\delta_m^j}|| + (\sum_{n\in I} ||(V_{j+1}\dots V_i h-2^j) 1_{\delta_n^{j+1}}||^2)^\frac12\\
\leq 3||(2^j-2^{j-1}) 1_{\delta_m^j}|| + 3(\sum_{n\in I} ||(V_{j+1}\dots V_i \underline{h} -2^j) 1_{\delta_n}||^2)^\frac12 = 3||(V_j\dots V_i \underline{h} -2^{j-1}) 1_{\delta_m^j}||.\endmultline$$ 

Finally, we have
$$||(V_8\dots V_i h -2^{7})^+1_{\delta_m^8}|| \leq 3 ||(V_8\dots V_i\underline{h} -2^{7})^+1_{\delta_m^8}||,\quad 0\leq m < 2^8,$$
which is more then we need.
\hfill $\square$

\enddemo

\proclaim{3.27. Lemma} For any triadic function $h$ we have
$$V h \leq ||(V_8\dots V_ih -2^7)^+||+ 2^8.$$

\endproclaim

\demo{Proof} Obviously
$$||V_0\dots V_ih|| \leq ||h_0|| +\dots + ||h_7|| + ||(V_8\dots V_i h -2^8)^+||\leq
 2+ (2+\dots + 2^7) + ||(V_8\dots V_ih -2^7||$$
with $h_0 =h\wedge 2^1$, $h_j =h\wedge 2^{j+1}- h\wedge 2^j$, $1\leq j <8$.
\hfill $\square$
\enddemo

\bigskip
\flushpar{\bf 3.28. Proof of Lemma 3.3.} By Lemmas 3.24-3.27, for some $U_j$ of type $j$, $j\geq 8$, the operations $W_j = V_jU_j$ are $j$-triadic and
$$\aligned ||V_8U_8\dots V_iU_i\underline{h_B} -2^8||&\geq (1-2^{-7})||V_8\dots V_i\underline{h_B} -2^8||-2^7\\
&\geq (1-2^7)(||V_8\dots V_i\underline{h_B} -2^7||-2^7)-2^{-7}\\
&\geq (1-2^{-7})\frac{1}{3} ||V_8\dots V_ih_B-2^7|| -2^7 +1 -2^{-7}\\
&\geq (1-2^{-7})\frac{1}{3} (V h_B-2^8) -2^7 +1 -2^{-7}.\endaligned$$
Thus the implication (21) is valid with $C = 644$, and Lemma 3.21 can be used.
\hfill $\square$

\bigskip

\flushpar{\bf 4. The construction of a discontinuous process for $V h_B =\infty$}. Up to now we have constructed orthogonal processes on {\it finite triadic sets} only (Lemma 3.3).

Now we need just one additional geometrical construction (cf. Definition 4.1). It showes in particular, that for any set $B$ satisfying (13) there exists a triadic set $\tilde B$ such that $h_{\tilde B}\geq h_B$, and the existence of an orthogonal a.e. discontinuous process on $B$ is equivalent to the existence of such a process on $\tilde B.$

Let us pass to the details. We start with the crucial
\definition{4.1. Definition} Let $A\subset [0, 1]$ be a set satisfying
$$\sharp A\cap [\alpha, 1] <\infty\quad\text{for any } \ \alpha >0.\tag22$$
The set $\tilde A$ {\it generated} by $A$ is defined by the formula
$$\tilde A = \{0, \frac{1}{3}, \frac{2}{3}, 1\}\cup\bigcup_{(i, n)\in I} \{n 3^{-2^i}, (n+1) 3^{-2^i}\},$$
where $I$ is a set of  pairs $(i, n)$, $i\geq 1$, $0\leq n <3^{2^i}$ for which there exist
$$t\in (n 3^{-2^i}, (n+1) 3^{-2^i}]\cap A$$
satisfying
$$\rho(t, A\setminus \{t\})\leq 3^{-2^{i-1}}.$$
\enddefinition

The following three lemmas can be obtained by easy and completely elementary considerations.

\proclaim{4.2. Lemma} Any set $\tilde A$ generated, by satisfying (22), $A$ is triadic.
\endproclaim

\proclaim{4.3. Lemma} For any set $\tilde B$ generated, by satisfying (13), $B$ we have $h_{\tilde B}\geq h_B.$
\endproclaim

\proclaim{4.4. Lemma} For $A, A_1$ satisfying  (22) and $A\subset A_1$ we have $\tilde A\subset\tilde A_1$. If additionally $\sharp A_1\setminus A <\infty$, then $\sharp \tilde A_1\setminus \tilde A <\infty$, in particular $\sharp \tilde A <\infty$ for $\sharp A <\infty.$
\endproclaim

Let $\rho(t, B) =\inf_{s\in B} |s-t|$. The following geometrical observation is particularly fruitful.

\proclaim{4.5. Lemma} For any (finite or countable) set $A\subset [0, 1]$ satisfying (22) and for its generated set $\tilde A$ we have
$$\sum_{t\in\tilde A} \rho(t, A) \leq 3,\tag23$$
$$\sum_{s\in A} \rho(s, \tilde A) \leq 1.\tag24$$
\endproclaim

\demo{Proof} Let $[\alpha, \beta] \cap A =\{\alpha, \beta\}$. Assume that
$$3^{-2^i} \leq\beta-\alpha <3^{-2^{i-1}},\quad i\geq 1.$$
Then
$$\sharp(\alpha,\beta)\cap \{n 3^{-2^j}; 0\leq j\leq i-1, 0\leq n \leq 3^{2^j}\}\leq 1$$
and
$$\sharp (\alpha, \beta)\cap\tilde A \cap \{n 3^{-2^j}; 0\leq n\leq 3^{2^j}\}\leq 2$$
for $j\geq i$. Thus
$$\multline\sum_{t\in\tilde A\cap [\alpha, \beta]} \rho(t, A) \leq \frac12 (\beta-\alpha) +2\cdot \sum_{j\geq i} 3^{-2^j}\\
< (\beta-\alpha) (\frac12 + 2\sum_{j\geq i} 3^{-2^j+2^i}) < 3 (\beta-\alpha).\endmultline$$

For $\beta -\alpha \geq 3^{-1}$ both the relation $\sum_{t\in\tilde A\cap [\alpha, \beta]} \rho(t, B) <3(\beta-\alpha)$  and the inequalities
$$\sum_{t\in\tilde A\cap [\beta_0, 1]}\rho(t, A) \leq 3 (1-\beta_0),$$
$$\sum_{t\in\tilde A\cap [0, \alpha_0]} \rho(t, A)\leq 3\alpha_0$$
for $\alpha_0 =\min A$, $\beta_0 =\max A$, are obvious. Hence (23).

On the other hand, $[\alpha, \beta] \cap A =\{\alpha,\beta\}$ implies $[\alpha,\beta]\cap \tilde A \ne \emptyset$ and
$$\sum_{s\in B\cap [\alpha,\beta]}\rho(s, \tilde A)\leq \beta -\alpha.$$
Moreover,
$$\sum_{s\in B\cap [0, \alpha_0]} \rho(s, \tilde A)\leq \alpha_0,$$
$$\sum_{s\in B\cap [\beta_0, 1]} \rho(s, \tilde A)\leq 1-\beta_0,$$
which yields (24).
\hfill $\square$
\enddemo

Let us remind the $X$ is an orthogonal process on $B$ satisfying (13) if $X: B\to L_2(\Bbb R)$, $X(0) =0$, $||X(t) - X(s)|| = |t-1|$ for $s,t\in B$. Then there obviously exists its extension $Y$ on $B\cup\tilde B$, being an orthogonal process on $B\cup\tilde B$. Moreover, the a.e. continuity of orthogonal processes on $B$ and their a.e. continuity on $\tilde B$ are equivalent. More precisely,

\proclaim{4.6. Lemma} For any orthogonal process $X$ on $B\cup\tilde B$ we have
$$X_{|B} \quad\text{is a.e. continuous if and only if} \
X_{|\tilde B} \quad\text{is a.e. continuous};$$
$$X_{|B} \quad\text{is a.e. discontinuous if and only if} \
X_{|\tilde B} \quad\text{is a.e. discontinuous.}$$
\endproclaim

\demo{Proof} For any function $\Delta: B\to L_2(\Bbb R)$ satisfying $\sum_{s\in B} ||\Delta(s)||^2 <0$, the a.e. continuity of $\Delta$ at $0$ is obvious.

Let $\phi: B\to\tilde B$, $\psi:\tilde B\to B$ be any functions satisfying
$$|s -\phi(s)| = \rho(s, \tilde B)\quad\text{for } s\in B,$$
$$|t-\psi(t)| =\rho(t, B)\quad\text{for } t\in\tilde B.$$
Then $\phi, \psi$ are non-decreasing and it is enough to look at $\Delta(s) = X(s) - X(\phi(s))$ on $B$ and then $\Delta(t) = X(t) -X(\psi(t))$ on $\tilde B$, cf. Lemma 4.5.
\hfill $\square$
\enddemo

As a corollary of Lemma 4.5 we also have 

\proclaim{4.7. Lemma} For any finite set $A\subset [0, 1]$, $ \alpha =\min A$, if $V h_{\tilde A} \geq C$ for some universal constant $C >0$, then there exists a process $X: A\to L_2[0, 1)$ with
$$X(\alpha) =0,\quad ||X(t) - X(s)||^2 = |t -s| \quad\text{for } s,t\in A,$$
satisfying $\langle X(t), 1_{[0, 1)]}\rangle =0$, $t\in A$, and
$$\lambda(\max_{t\in A} |X(t)| > 1) > \frac{1}{6}.$$
\endproclaim

\demo{Proof} There exists an orthogonal process $X_1: \tilde A \to L_2(\Bbb R)$ satisfying
$$\lambda([0, 1)\cap (\max_{s\in \tilde A} X_1(s)^2 >75)) >\frac12$$
if only $C = 5\sqrt{3} c$ with $c$ given by Lemma 3.3. Taking $X = u X_1$ with a suitable unitary operator $u: L_2(\Bbb R)\to L_2 [0, 1)$, we have also
$$\lambda(\max_{s\in\tilde A} X(s)^2 >75) >\frac{1}{3},$$
for some orthogonal process with values in $L_2[0, 1)$. One can also assume that $X$ has been extended to an orthogonal process $X: \tilde A\cup A\to L_2[0, 1)$, and that $\langle X(t), 1_{[0, 1)}\rangle =0$, by standard tricks.

Let, analogously to 4.6, $\rho(s, A) = |s -\phi(s)|$, $s\in \tilde A$, $\phi: \tilde A\to A$.
Schwartz and Tshebyshev's inequalities lead to the following rather obvious estimates. For any $s\in\tilde A$,
$$X(s)^2\leq 3 [(X(\phi(s))- X(\alpha))^2 + (X(\phi(s)) - X(s))^2 + x(\alpha)^2]$$
and
$$\max_{t\in A} (X(t) - X(\alpha))^2\geq \max_{s\in\tilde A} (X(\phi(s)) - X(\alpha))^2\geq \frac{1}{3} \max_{s\in \tilde A} X(s)^2 - Y$$
for
$$Y = X(\alpha)^2 +\sum_{s'\in \tilde A} (X(\phi(s')) - X(s'))^2.$$
By Lemma 4.4, we have
$$\int_{[0, 1)} Y\leq 4$$
and
$$\lambda(Y \geq 24) <\frac{1}{6}$$
and, by assumptions on $X(s)$, $s\in\tilde A$,
$$\lambda(\max_{t\in A} (X(t) - X(\alpha))^2 >1) > \frac{1}{3}-\frac{1}{6}.$$
\vskip-12pt
\hfill $\square$
\enddemo

\baselineskip18pt

From now on let $B$ be a fixed set satisfying (13), with$V h_B =\infty.$

\proclaim{4.8. Lemma} For any $\beta\in B$ there exists $0 <\alpha <\beta$, $\alpha\in B$, satisfying $V h_{B\cap [\alpha, \beta]} \geq C$, with $C$ being a given constant.
\endproclaim

\demo{Proof} By Lemmas 4.3, 4.4,  $V_{\tilde B} = \infty$, $\sharp\tilde B\setminus \widetilde{B\cap [0, \beta]} <\infty$ thus $(h_{\tilde B}-2^i)^+ = (h_{\widetilde{B\cap [0, \beta]}} -2^i)^+$ for some $i\geq 1$ and $V h_{\widetilde{B\cap [0,\beta]}} = \infty$. In particular
$$||V_0\dots V_k h_{\widetilde{B\cap [0,\beta]}} || > C +1$$
for some $k\geq 1$, and by the Lebesgue monotone convergence theorem  (used $k+1$ times), we also have 
$$||V_0\dots V_k h_{\widetilde{B\cap [\alpha, \beta]}}|| >C$$
for sufficiently small $\alpha >0.$

\enddemo

\bigskip
\flushpar{\bf 4.9. The construction of a discontinuous process}. Let us recall that $B$ is a fixed set satisfying (13) and $V h_B=\infty$. By Lemmas 4.7, 4.8, there exist numbers $1 =\alpha_1 >\alpha_2 >\dots >0$ satisfying
$$V h_{\widetilde{B\cap [\alpha_{s+1}, \alpha_s]}} >C, \quad s\geq 1,$$
and
$$\lambda(\max_{t\subset B\cap [\alpha_{s+1}, \alpha_s]} |X_s(t)| >1) >\frac{1}{6}$$
for some processes $X_s:\! B\cap [\alpha_{s+1}, \alpha_s]\to L_2[0, 1)$, $X_s(\alpha_{s+1}) =0$, $||X_s(t) - X_s(t_1)||^2 = |t-t_1|$, $\langle X(t), 1_{[0, 1)} \rangle =0.$
 Taking the space $([0, 1)^\Bbb N, \lambda^{\otimes\Bbb N})$, isomorphic to $([0, 1),\lambda)$, and denoting
$$\omega = (\omega_1, \omega_2,\dots)\quad\text{for } \omega\in [0, 1)^\Bbb N,$$
it is enough put
$$X(t)(\omega) = X_s(t)(\omega_s) +\sum_{s' >s} X_{s'}(\alpha_{s'})(\omega_{s'})$$
for $t\in (\alpha_{s+1}, \alpha_s]$. By Borel-Cantelli Lemma, the Cauchy condition for $X(t)$, $t\to 0$, $t\in B$, fails on a set $Z$ of measure $\lambda^{\otimes\Bbb N}(Z) =1.$
The existence of an a.e. discontinues process with values in $L_2[0, 1)$ is equivalent to the  existence of such a process with values in $L_2(\Bbb R).$

\bigskip

\flushpar{\bf 5. The proof of continuity of any process  for $V h_B<\infty$}.
Assume now that $V h_B <\infty$ for some set $B$ satisfying (13). 
We show that any orthogonal process $X$ on $B$ is a.e. continuous. By Lemmas 4.2, 4.6,
it is enough to show that each such process on $\tilde B$ is continuous a.e. for the triadic generated set $\tilde B$.  Thus the proof of the continuity of an orthogonal process $X$ on $B$ splits into two main parts.

First we show that each such process is a.e. continuous  if $B$ is a {\it triadic} set with $V h_B <\infty$ (cf. Lemma 5.7).
The proof is pretty simple and based on classical ideas (due to Plancherel and Tandori). 

More tedious estimates for information functions are needed in the second part, when we show that $V h_B <\infty$ implies $V h_{\tilde B} <\infty$ for any set $B$ satisfying (13) and its generated set $\tilde B$ (Lemmas 5.10, 5.12).

In the first part a crucial role is played by the classical Plancherel idea of diadic partitions of a given sequence of vectors. It gives (see [2])

\proclaim{5.1. Lemma} For any othogonal vectors $Y_1,\dots Y_N$ in $L_2$, we have
$$||M||^2 \leq k^2 \sum_{1\leq n\leq N} ||X||^2$$
for
$$M =\max_{1\leq n\leq N} |Y_1 +\dots + Y_N|,$$
$$k \geq \log_2 N+ 1.$$
\endproclaim

In particular, if an orthogonal process $X$ is defined on the whole set
$$A^j_m = \{n 3^{-2^{j+1}};  0\leq n < 3^{2^{j+1}}\}\cap \delta_m^j\tag25$$
for a fixed $j\geq 0$, $0\leq m <3^{2^j}$, then
$$||\max_{t\in A_m^j}|X(t) - X(m 3^{-2^j})| \ ||\leq 3\cdot 2^j\sqrt{\lambda(\delta_m^j)}.\tag26$$

Let us fix  a finite triadic set $B$ and an orthogonal process $X$ on $B$.

\medskip
\flushpar{\bf 5.2. Notation.} For any $j\geq 1$, $0\leq m < 3^{2^j}$, denote
$$M_m^j = \max_{t\in\delta_n^j\cap B} |X(t) - X(n 3^{-2^j})|$$
if $(n 3^{-2^j}, (n +1)3^{-2^j})\cap B\ne \emptyset$, and $M_m^j =0$ otherwise.

\medskip
\flushpar{\bf 5.3. Notation.} For any triadic function $h$ we put
$$\bar V_j h = (h\wedge 2^j) +2^j 1_{(h\geq 2^j)} + ||(h -2^{j+1})^+||_j,$$
which is obviously another triadic function.

The operations $\bar V_j$ are now  useful because of the following crucial estimate.

\medskip

\proclaim{5.4. Lemma} Assume that for some triadic function $h$ and for $j\geq 0$ we have $h\wedge 2^{j+1} = h_B\wedge 2^{j+1}$ and, for $M_n^{j+1}$ defined in 5.2,
$$||M_n^{j+1}|| \leq 3||(h -2^{j+1})^+ 1_{\delta_n^{j+1}}||$$
for any $0\leq n < 3^{2^{j+1}}$. Then
$$||M_m^j|| \leq 3||(\bar V_j h-2^j)^+ 1_{\delta_m^j}||$$
for $0\leq m < 3^{2^j}$, and obviously  $\bar V_j h\wedge 2^j = h_B \wedge 2^j.$
\endproclaim

\demo{Proof} Let us fix $\delta_m^j$ such that Int $\delta_m^j\cap B\ne\emptyset$. To unify our considerations, let us assume that $X$ is extended is such  a way that $A_m^j\subset B$ for the set given by (25). Using (26), we get
$$||M_m^j||\leq ||\max_{t\in A_m^j} |X(t) - X(m 3^{-2^j})| \ || + ||\max_{\smallmatrix 0\leq n <3^{2^{j+1}}\\ \delta_n^{j+1}\subset\delta_m^j\endsmallmatrix} M_n^{j+1}||,$$
$$||\max_{t\in A_m^j} |X(t) - X(m 3^{-2^j})| \ ||\leq 3 ||2^j 1_{\delta_m^j}||,$$
$$\multline ||\max_{\smallmatrix 0\leq n <3^{2^{j+1}}\\ \delta_n^{j+1}\subset\delta_m^j\endsmallmatrix}  M_n^{j+1}||^2\leq \int_\Bbb R \sum_{\smallmatrix 0\leq n <3^{2^{j+1}}\\ \delta_n^{j+1}\subset\delta_m^j\endsmallmatrix} (M_n^{j+1})^2\\
\leq 9 \sum_{\smallmatrix 0\leq n <3^{2^{j+1}}\\ \delta_n^{j+1}\subset\delta_m^j\endsmallmatrix} ||(h -2^{j+1})^+ 1_{\delta_n^{j+1}}||^2 = 9||(h- 2^{j+1})^+ 1_{\delta_m^j}||^2.\endmultline$$
By the definition of conditional norm, the equality
$$||2^j 1_{\delta_m^j}|| + ||(h-2^{j+1})^+ 1_{\delta_m^j}|| = ||(\bar V h -2^j)^+ 1_{\delta_m^j}||$$
is obvious; in fact, in the case Int $\delta_m^j\cap B\ne\emptyset$ we have $\delta_m^j\subset (h_B\geq 2^j) = (h\geq 2^j)$.

In the oposite case there is nothing to prove.\hfill $\square$

\enddemo

By  an obvious backward induction we have

\proclaim{5.5. Corollary} If a finite triadic set $B$ satisfies $h_B\leq 2^{i+1}$, $i\geq 0$, then
$$||M_n^j|| \leq 3 ||(\bar V_j\dots  \bar V_i h_B-2^j)^+ 1_{\delta_n^j}||$$
for $0\leq j\leq i$, $0\leq n <3^{2^j}$, for any orthogonal process $X$ on $B$ and for $M_n^j$ given by {\rm 5.2.}

\endproclaim

In fact we can use operations $V_j$ instead of $\bar V_j$ because of the following trick.

\proclaim{5.6. Lemma} For any triadic function $h$ satisfying $h\leq 2^{i+1}$, we have
$$(\bar V_j\dots \bar V_i h-2^j)^+ \leq (2 V_j\dots V_i h-2^j)^+,\quad 0\leq j < i.$$
\endproclaim

\demo{Proof} We leave to the  reader the fairly straightforward backward induction showing that
$$(\bar V_j\dots \bar V_i h-2^j)^+ = (2 V_j\dots V_i \frac12 \bar h-2^j)^+$$
with $\bar h = 2^{j+1}$ for $2^j\leq h < 2^{j+1}$, $0\leq j <i$, $\bar h=2^{i+1}$ for $2^i\leq h \leq 2^{i+1}$. Obviously $\frac12 \bar h \leq h.$ \hfill $\square$ 
\enddemo

The continuity of orthogonal processes on {\it any} triadic set (the goal of the first part of this section) can be obtained  now by a simple formal reasoning.

\proclaim{5.7. Lemma } For any triadic set $B$ with $V h_B <\infty$, any orthogonal process $X(t)$, $t\in B$, is continuous.
\endproclaim

\demo{Proof} We have
$$ ||\lim_{i\to\infty} V_j\dots V_i h_B 1_{\delta_0^j}||\to 0\quad
\text{for } \quad j\to \infty.$$
Thus, for suitably chosen
$j(1) < j(2) <\dots,$
we have, in particular,
$$\sum_{s\geq 1} ||\lim_{i\to\infty} (2 V_{j(s)}\dots V_i h_B-2^{j(s)})^+|| <\infty.\tag27$$
For suitably chosen $i(s) >j(s)$, $s\geq 1$, denoting $B_s =\{n 3^{-2^{i(s)+1}}; 0\leq n\leq 3^{2^{i(s)+1}}\}$, we have 
$$ B\setminus \delta_0^{j(s+1)} \in B_s\quad\text{(cf. Definition 3.1)},$$
and
$$M(s):= \max_{t\in B\cap(\delta_0^{j(s)}\setminus \delta_0^{j(s+1)})} |X(t)| \leq M_0^{j(s)}$$
for $M_0^{j(s)}$ given in Definition 5.2 in which, instead of $X$, we take the process $X_{|B\cap B_s}$ restricted to a finite triadic set $B\cap B_s$. Obviously, $h_{B\cap B_s}\leq 2^{i(s)+1}$ and, using Lemma 5.6 and Corollary 5.5, we get
$$\aligned ||M(s)||\leq ||M_0^{j(s)}|| \leq & 3 ||(2 V_{j(s)}\dots V_{i(s)} h_{B\cap B_s}-2^{j(s)})^+||\\
= & 3||\lim_{i\to\infty} (2 V_{j(s)}\dots V_i h_{B\cap B_s}-2^{j(s)})^+||\\
\leq & 3||\lim_{i\to\infty} (2 V_{j(s)}\dots V_i h_B -2^{j(s)})^+||.\endaligned$$
By (27), this finishes the proof.
\hfill $\square$

\enddemo

Let $\tilde B$ be a triadic set generated by the given set $B$ satisfying (13). Thus $\tilde B$ can be used in Lemma 5.7. By Lemma 4.6, continuity of orthogonal processes on $B$ and on $\tilde B$ are equivalent. Thus it is enough to show that the inequality $V h_{\tilde B} <\infty$ is implied by $V h_B <\infty$. It can be done in a number of elementary ways. Precise calculations are fairly  tedious.

\proclaim{5.8. Lemma} Let $V g = g\wedge 4 + ||(g -4)^+||$  act in a real $L_2$-space for a probability measure $P$.
Let moreover $g\geq g_1\geq 2$ and, for $A = (g\geq 8)$,
$$||(g -2) 1_A||\leq 14||(g_1 -2) 1_A||,\quad g_1\geq 4 \ \text{on } A,\tag28$$
for  some given elements $g, g_1\in L_2$. Then
$$||V g-1||\leq 14||V g_1 -1||$$
and, obviously,
$$V g\geq V g_1\geq 2.$$
\endproclaim

\demo{Proof}  Put $x= ||(g_1-4) 1_A||$, $\Delta = P(A)$, $\Delta' = P(A^c)$. Then, by (28),
$$\aligned ||V g-1|| &\leq 3+ ||(g -4)^+||\leq 3 +4 + ||(g -8) 1_A||\\
& \leq 7+ 14||(g_1 -2) 1_A||\\
& \leq 7 + 14(x +2\sqrt{\Delta}) = 14(x +2\sqrt{\Delta} +\frac12).\endaligned$$
On the other hand
$||(g_1 -4)^+||\geq  x$, $V g_1\geq (4+x) 1_A+ (2 +x) 1_{A^c}$ and
$$\multline ||V g_1 -1|| \geq \sqrt{(x+3)^2\Delta + (x+1)^2\Delta'}\\
= ((x+2\sqrt{\Delta}+\frac12)^2 + (\sqrt{x} - 2\sqrt{x \Delta})^2 + (\frac12 - 2\sqrt{\Delta})^2 + \frac12 + 4x\Delta)^{\frac12}\geq x + 2\sqrt{\Delta} +\frac12.\endmultline$$
\vskip-12pt
\hfill $\square$
\enddemo

\baselineskip18pt

\medskip
\flushpar{\bf 5.9. Notation.} For any number $a\geq 1$ we write
$$\underset{=}\to{a} = 2^{j-1}\quad\text{for } 2^j\leq a < 2^{j+1}, \ j\geq 0.$$

\proclaim{5.10. Lemma} For any triadic function $h$ satisfying $h\leq 2^{i+1}$ we have
$$||V_0\dots V_i h||\leq 14||V_0\dots V_i \underset{=}\to{h}||.$$

\endproclaim

\demo{Proof} We show that $||(V_j\dots V_i h-2^{j-2}) 1_{\delta_m^j}||\leq 14  ||(V_j\dots V_i \underset{=}\to{h} - 2^{j-2})1_{\delta_m^j}||$ by backward induction on $j$, $0\leq j\leq i$.  For any interval $\delta_n^i = (n 3^{-2^i}, (n+1) 3^{-2^i}]\subset (h\geq 2^i)$, we have
$$\underset{=}\to{h} \geq 2^{i-1}\quad\text{on}\quad \delta_n^i,\tag29$$
$$\aligned ||(V_i h-2^{i-2}) 1_{\delta_n^i}||\leq & 7 ||2^{i-2} 1_{\delta_n^i}||\\
\leq &7 ||(V_i\underset{=}\to{h} -2^{i-2}) 1_{\delta_n^i}||\\
\leq & 14||V_i \underset{=}\to{h} -2^{i-2}) 1_{\delta_n^i}||.\endaligned\tag30$$

Let us fix $j$, $0\leq j <i$ and suppose that
$$V_{j+1}\dots V_i \underset{=}\to{h} \geq 2^j$$
on any interval $\delta_n^{i+1}\subset (V_{j+1}\dots V_i h\geq 2^{j+1}),$
and that
$$||(V_{j+1}\dots V_i h-2^{j-1})1_{\delta_n^{j+1}}||
\leq 14||(V_{j+1}\dots V_i\underset{=}\to{h} -2^{j-1}) 1_{\delta_n^{j+1}}||$$
(thus we assume (29), (30) with $j+1$, $V_{j+1}\dots V_ih$, $V_{j+1}\dots V_i \underset{=}\to{h}$ instead of $i, h, \underset{=}\to{h}$).

By the definition of a triadic function, for
$$\delta_m^j \subset(h\geq 2^j),\quad A =\delta_m^j\cap (h\geq 2^{j+1}),$$
we have
$$A =\bigcup_{1\leq k\leq K} \delta_{n(k)}^{j+1}$$
for some $0\leq n(k) < 3^{2^{j+1}}$. Thus all assumptions of Lemma 5.8 are satisfied for $L_2 = L_2(\delta_m^j)$, with the normalized Lebesgue measure,
$$g = 2^{-j+2} V_{j+1}\dots V_i h\quad\text{on } \ \delta_m^j,$$
$$g_1 = 2^{-j+2} V_{j+1}\dots V_i \underset{=}\to{h}\quad\text{on } \ \delta_m^j.$$
Thus
$$V_j\dots V_i \underset{=}\to{h} \geq 2^{j-1}\quad\text{on } \ \delta_m^j,$$
$$||(V_j\dots V_i h-2^{j-2})1_{\delta_m^j}||\leq 14||(V_j\dots V_i\underset{=}\to{h} -2^{j-2})1_{\delta_m^j}||.$$

By backward induction the last inequality is valid for any $0\leq j\leq i.$ Taking $j=0$ we obtain more then needed. \hfill $\square$
\enddemo

The proof of Theorem 1.2 is now reduced to the implication $V h_B <\infty\implies V\underset{=\!=}\to{h_{\tilde B}} <\infty$. Unexpectedly, for arguments $h_B$, $\underset{=\!=}\to{h_{\tilde B}} $,
being not necessarily triadic functions,  the comparison of values $V_i h_B$, $V_i \underset{=\!=}\to{h_{\tilde B}} $ of our crucial operation $V_i$ is complicated. To compress the calculations, we use the following notation (cf. 1.4).

\medskip
\flushpar{\bf 5.11. Notation.} For any function $f\geq 1$, we write
$$f_{j\downarrow} = f\wedge 2^j$$
and
$$f_j = f_{j+1\downarrow} -f_{j\downarrow}, \quad f_{j\uparrow} = f -f_{j\downarrow},$$
for $j\geq 0$, in particular
$$f = f_{j\downarrow} + f_j + f_{j+1\uparrow}$$
for any real function $f.$

\proclaim{5.12. Lemma} For any Borel functions $1\leq g\leq h$ on $(0, 1]$ and for
$$A =\bigcup_{\smallmatrix 0\leq n < 3^{2^j}\\ \delta_n^j\cap(h\geq 2^j)\ne\emptyset\endsmallmatrix} \delta_n^j\tag31$$
with $j\geq 1$, we have
$$ ||(V_j h)_{j-1\uparrow} - (V_j g)_{j-1\uparrow}||
 \leq ||h_{j\uparrow} -g_{j\uparrow}|| + ||(h_{j\downarrow} -g_{j\downarrow}) 1_A|| + ||h_{j-1}- g_{j-1}||.$$

\endproclaim

\demo{Proof} The required calculations are natural:
$$\aligned (V_jh)_{j-1\uparrow} & - (V_j g)_{j-1\uparrow}\\
& = (h_{j-1} - g_{j-1})1_A c + ((V_j h)_{j-1\uparrow} - (V_j g)_{j-1\uparrow}) 1_A\\
& \leq (h_{j-1} - g_{j-1}) 1_A c + (V_j h - V_jg) 1_A\\
&\leq (h_{j-1}- g_{j-1}) 1_A c + ||h_{j\uparrow} - g_{j\uparrow}||_j + (h_{j\downarrow} -g_{j\downarrow})1_A\\
& \leq (h_{j-1}-g_{j-1}) + ||h_{j\uparrow} -g_{j\uparrow}||_j + (h_{j\downarrow} -g_{j\downarrow}) 1_A.\endaligned$$
\enddemo

\proclaim{5.13. Lemma} For any set $B$ satisfying (13) and for $0\leq j <i$ let us put
$$g = V_{j+1}\dots V_i h_B,$$
$$h = V_{j+1} \dots V_i (h_B\vee \underset{=\!=}\to{h_{\tilde B}} )$$
and let the set $A$ be given by (31). Then
$$||h_{j-1} - g_{j-1}|| \leq 2^j\cdot 3^{-2^{j-1}},\tag32$$
$$||(h_{j\downarrow} -g_{j\downarrow}) 1_A||\leq 2^{j+1} 3^{-2^{j-2}}.\tag33$$
\endproclaim

\demo{Proof} The inequality (32) is an immediate consequence of our Definition 4.1. Let us fix an interval $(\alpha, \beta]$ satisfying $[\alpha, \beta]\cap B = \{\alpha, \beta\}$. Observe that, by 5.9,
$$(\alpha, \beta]\cap (\underset{=\!=}\to{h_{\tilde B}}  >2^{j-1}) = (\alpha, \beta]\cap (h_{\tilde B} > 2^{j+1})\subset \delta_k^{j+1}\cup\delta_l^{j+1}$$
for some $0\leq k\leq l < 3^{2^{j+1}}$, and, because of the structure of $V_{j+1},\dots, V_i$,
$$(\alpha, \beta]\cap (h > 2^{j-1}) \subset \delta_k^{j+1}\cup \delta_l^{j+1}.$$
Thus $\lambda((\alpha, \beta]\cap [h > 2^{j-1})) < 2\cdot 3^{-2^{j+1}}.$

On the other hand
$$\multline \sharp \{(\alpha, \beta];\quad [\alpha, \beta]\cap B =\{\alpha, \beta\}, (g < 2^j)\cap (\alpha, \beta]\ne \emptyset\}\\
\leq \sharp \{((\alpha, \beta]; \quad [\alpha, \beta]\cap B =\{\alpha, \beta\}, (h_B <2^j)\cap (\alpha, \beta]\ne\emptyset\}\leq 1/3^{-2^j},\endmultline$$
as $(h_B)_{j+1\downarrow} \leq g_{j+1\downarrow}$. Inequality (32) is given by a natural estimate
$$||h_{j-1}-g_{j-1}||\leq 2^{j-1} [\bigcup_{\smallmatrix (\alpha, \beta]\\ [\alpha, \beta] \cap B =\{\alpha, \beta\}\endsmallmatrix} \lambda((\alpha, \beta]\cap (g <2^j)\cap (h >2^{j-1}))]^{\frac12}.$$

To prove (33) let us fix $\delta_m^j\subset A$, $0\leq m < 3^{2^j}$. Supposing $\delta_m^j\cap B =\emptyset$, we have $\delta_m^j\cap \tilde B =\emptyset$ and $\delta_m^j\subset A^c$. A careful analysis is needed only when
$$k =\min \{n; \delta_n^{j+1}\cap\delta_m^j\cap B\ne \emptyset\},$$
$$l= \max \{n; \delta_n^{j+1}\cap \delta_m^j\cap B\ne \emptyset\}$$
are defined. The case $k= l$ is obviously possible. Then 
$$h_B >2^j \quad\text{on } \bigcup_{k <n<l} \delta_n^{j+1}$$
and, by the structure of $V_{j+1},\dots, V_i$,
$$h\geq g >2^j\quad\text{on } \bigcup_{k <n < l} \delta_n^{j+1}.$$
Moreover 
$$(B\cup\tilde B)\cap  \delta_m^j\cap \bigcup_{n <k \ \text{or} \ n >l} \delta_n^{j+1}=\emptyset$$
and $h_B$, $h_{\tilde B}$ are constant on the set $\delta_m^j\cap \bigcup_{n <k}\delta_n^{j+1}$ as well as on the set $\delta_m^j\cap \bigcup_{n >l} \delta_n^{j+1}$. It can be easily verified that
$$\underset{=\!=}\to{h_{\tilde B_1}} \leq 2^{j-1}\quad\text{on } \delta_m^j\cap \bigcup_{n <k \ \text{or } n >l} \delta_n^{j+1},$$
and obviously $g= h_B$, $h =h_B\vee \underset{=\!=}\to{h_{\tilde B_1}}$, on $\delta_m^j\cap \bigcup_{n < k \ \text{or} \ n >l} \delta_n^{j+1}$. We say that $m$ is {\it special} if
$$(g < 2^{j-2})\cap \delta_m^j \cap \bigcup_{n <k \ \text{or } n >l} \delta_n^{j+1} \ne \emptyset.$$
Then, by natural estimation of the number of "large" intervals $(\alpha, \beta]$ defined by 
 $[\alpha, \beta]\cap B =\{\alpha, \beta\}$, $\beta -\alpha >3^{-2^{j-1}}$, we have
$$\sharp\{0\leq m <3^{2^j}; m  \ \text{is special}\} \leq 2 \sharp \ \text{(set of "large" intervals $(\alpha, \beta])$},$$
and
$$\multline||h_{j\downarrow} - g_{j\downarrow}||^2 \leq  2^{2j} \cdot 3^{-2^{j+1}} \cdot 2\cdot\sharp \{0,\dots, 3^{2^j}-1\}\\
+ 2^{2j}\cdot 3^{-2^j}\cdot\sharp\{0\leq m < 3^{2^j}; m \  \text{is special}\}\leq 2\cdot 2^{2j}(3^{-2^j} +3^{-2^j}\cdot 3^{2^{j-1}}).\endmultline$$
Hence (33).\hfill $\square$
\enddemo

\medskip

\flushpar{\bf 5.13. Proof of Theorem 1.2}. For any set $B$ satisfying (13) with $h_B\leq K < \infty$ we have $V h_{B_1}\leq K$ for any finite subset $B_1\subset B$, satisfying (13). 
Obviously $h_{B_1}, \underset{=\!=}\to{h_{\tilde B_1}}\leq 2^i$ for suitably large $i \geq 1$, and Lemmas 5.11, 5.12 give inequalities
$$\aligned  ||V_0\dots & V_i \underset{=\!=}\to{h_{\tilde B_1}} -1||\\
& \leq ||V_0\dots V_i(\underset{=\!=}\to{h_{\tilde B_1}} \vee h_{B_1}) -1||\\
&\leq ||V_0\dots V_i h_{B_1} -1|| +\sum_{j\geq 0} (2^j\cdot 3^{-2^{j-1}} + 2^{j+1}\cdot 3^{-2^{j-2}})\\
& \leq K +\sum_{j\geq 0} (2^j\cdot 3^{-2^j}+2^{j+1}\cdot 3^{-2^{j-2}}) =: L,\endaligned$$
and, by Lemma 5.10,
$$||V_0\dots V_i h_{\tilde B_1} || \leq 14(L +1).$$
Then $V h_{\tilde B} < 14(L +1)$ by an obvious application of the Lebesgue monotone convergence theorem, and Lemmas 5.7, 4.6 give a.e. convergence of any orthogonal process on $B$.

Together with 4.9, it completes the proof of 1.2.

\bigskip
\flushpar{\bf 6. Applications to the continuity of orthogonal processes}

The definition of an orthogonal process $X$ on $B$ can be used for any set $B$, $\{0, 1\}\subset B\subset [0, 1]$. Taking $X(t) =\lim_{n\to\infty} X(t_n)$ in $L_2$, with some $t_n\to t$, $t_n\in B$, for any $t\in \bar B$, one can extend $X$ to an orthogonal process defined on the closure $\bar B$ of $B$. 

It is natural to assume in what  follows that the set of time moments $B$ is closed and $\{0, 1\}\subset B\subset [0, 1]$. In such a general case it is natural to accept infinite values of information function of the 'partition' given by $B$. Thus, we define
$$\aligned H_B(t) = &\infty \quad\text{for}\quad t\in B,\\
H_B(t) = & -\log_3(\beta -\alpha)\quad\text{for}\quad t\in (\alpha, \beta), [\alpha, \beta]\cap B =\{\alpha, \beta\}.\endaligned\tag34$$

Obviously, $H_B = h_B$ a.e., if only $B$ satisfies (13). 

As usually, $\tilde X$ is a {\it version} of the process $X$ if $\tilde X(t)$ represents an element $X(t)\in L_2$ and $\tilde X(t)$ is a function, for any $t\in B$. Allowing possibly uncountable sets $B$, we use the following definition:

$X: B\to L_2$ is a.e. continuous at $t\in B$ if there exists a version of $X$ a.e. continuous at     
 $t.$

$X: B\to L_2$ is a.e. discontinuous at a point $t\in B$ if any version of $X$ is a.e. discontinuous at $t.$

As a rather easy and formal corollary we obtain now the main result of this section.

\proclaim{6.1. Theorem} For any $t\in B = \bar B$, $\{0, 1\}\subset B\subset [0, 1]$:

{\rm A.} Any orthogonal process on $B$ is a.e. continues at $t$ if and only if $V(H_B 1_U) <~\infty$ for some neighbourhood $U$ of $t.$

{\rm B.} There exists an orthogonal process on $B$ a.e. discontinues at $t$ if and only if $V(H_B 1_U) =\infty$ for any neighbourhood $U$ of $t.$

\endproclaim

\definition{6.2. Definition} Let us say that a Lebesgue measure preserving, one to one mapping $S^j: (0, 1]\to (0, 1]$ is of $j$-type, $j\geq -1$, if
$S^j$ restricted to any $\delta_n^{j+1}$ is a shift by a multiple of $3^{-2^{j+1}}$, i.e.
$$S^j(t) = t+ k\cdot 3^{-2^{j+1}}\quad\text{for } \ t\in  \delta_n^{j+1}$$
with an integer $k$, depending on $0\leq n <3^{2^{j+1}}.$

Moreover, we say that some $j$-type mapping $S_m^j$, $j\geq 0$, is of $j, m$-type, $0\leq m < 3^{2^j}$, if  $S_m^j$ is the identity on $(0, 1]\setminus \delta_m^j$. Finally, for $t\in [0, 1]\setminus \{n 3^{-2^j}; j\geq 0, 0\leq n \leq 3^{2^j}\}$,
 the mapping $S_t: [0, 1]\to (0, 1]$ is a $t$-operation if
$$S_t=\lim_{i\to\infty} S^{-1} S_{m_0}^0\circ \dots \circ S_{m_i}^i$$
for some mapping $S^{-1}$ of $(-1)$-type and mappings $S_{m_j}^j$ being of $j, m_j$-types, respectively, with $m_j$ satisfying $\delta_{m_0}^0\supset \delta_{m_1}^1\supset\dots,  \bigcap_{j\geq 0} \delta_{m_j}^j =\{t\}.$ 
\enddefinition

As $\sum_{j\geq 0} \sup_{t\in (0, 1]} |S_{m_j}^j t-t| \leq \sum_{j\geq 0} 3^{-2^j} < \infty$, the limit $S_t$ is a well defined  Lebesgue measure preserving mapping, one to one on $(0, 1]\setminus \{t\}.$

\proclaim{6.3. Lemma} For any set $B =\bar B$, $\{0, 1\}\subset B\subset [0, 1]$ and $t$-operation $S$ we have $V H_B = V H_{\overline{S(B)}}$.

\endproclaim

\demo{Proof} Note that $||V_0\dots V_i H_{\bar B_1}|| = ||V_0\dots V_i H_{\overline{S_m^j B_1}}||$ for any $\{0, 1\}\subset B_1\subset [0, 1]$ and $S_m^j$ being of $j, m$-type. Thus, by the monotonic passage  to the limits,
$$||V_0\dots V_i H_B|| = ||V_0\dots V_i H_{\overline{S^{-1}S_{m_0}^0\dots S_{m_i}^i(B)}}||,$$
$$||V_0\dots V_i H_B|| = ||V_0\dots V_i H_{\overline{S(B)}}||,$$
$$V H_B = V H_{\overline{S(B)}}.$$

\enddemo

\proclaim{6.4. Lemma} For any $t\in (0, 1]\setminus \{n 3^{-2^j}; j\geq 0, 0\leq n\leq 3^{2^j}\}$ there exist $m_j, j\geq 0$, and $S_{m_j}^j$ of $j, m_j$-type, $j\geq 0$, and $S^{-1}$ of $(-1)$-type such that $S_t(t) =0$ for the $t$-operation $S$ given by Definition 6.2, and $S_t$ is continuous at $t.$

\endproclaim

\demo{Proof} Define $m_j$  by the condition $t\in \delta_{m_j}^j$, $j\geq 0$, and put
$$ S_{m_j}^j(\delta_{n_1}^{j+1}) = \delta_{n_3}^{j+1},$$
$$S_{m_j}^j(\delta_{n_3}^{j+1}) = \delta_{n_1}^{j+1}$$
for $n_1, n_2, n_3$ defined by
$$\delta_{m_j}^j =\bigcup_{n_1\leq n\leq n_2} \delta_n^{j+1},\quad t\in \delta_{n_3}^{j+1}.$$
$S^{-1}$ can be defined by a similar trick. \hfill $\square$

\enddemo

By Lemmas 6.3, 6.4, the investigation of continuity of an orthogonal process in any point $t$ can be reduced to the case $t=0.$

\proclaim{6.5. Lemma} Let $S_t$ be defined as in 6.2 for some $t\in [0, 1]\setminus \{n 3^{-2^j}; j\geq 0, 0\leq n\leq 3^{2^j}\}$, and let $S_t(0) =0$.

{\rm A.} Any orthogonal process on $B$ is a.e. continuous at $t$ if and only if each orthogonal process on $S_t(B)$ is continuous at $0.$

{\rm B.} There exists an orthogonal process on $B$ a.e. discontinuous at $t$ if and only if there exists an orthogonal process on $S(B)$ a.e. discontinuous at $0.$

\endproclaim

\demo{Proof} Observe that for any orthogonal process $X$ on $B$ and any mapping $S_m^j$ of $j, m$-type, one can put
$$Y(S_m^j(s)) = X(s) + R(s),\quad s\in B,$$
being an orthogonal process on $S_m^j(B)$ for some $R: (0, 1]\to L_2(\Bbb R)$,
$$R \ \text{being constant on each } \delta_n^{j+1},$$
$$||R(t)||^2 \leq 3^{-2^j}\quad \text{for any } t\in (0, 1].$$

Consequently, for any process $X$ orthogonal on $B$ there exists a process $R: (0, 1]\to L_2(\Bbb R)$ such that
$$Y(S_t(s)) = X(s) + R(s),\quad s\in B,$$
is orthogonal on $S_t(B)$ and
$$R \ \text{is a.e. continuous at } t.$$

Similarly, for any process $Y$ orthogonal on $S_t(B)$ there exists a process $P: [0, 1]\to L_2(\Bbb R)$ for which
$$X(s) = Y(S_t(s)) + P(S_t(s)),\quad s\in B,$$
is orthogonal on $B$ and
$$P \ \text{is  a.e. continuous in } 0.$$
The construction of the processes $R$ and $P$ gives both equivalences $A$ and $B.$

\enddemo

\proclaim{6.6. Lemma} {\rm A.} Any orthogonal process on $B$ has a version a.e. continuous at $t =0$ if and only if $V(H_B 1_{[0, \frac{1}{k}]}) <\infty$ for some $k\geq 1.$

{\rm B.} There exists an orthogonal process on $B$ with all versions a.e. discontinuous at $t =0$ if and only if $V(H_B 1_{[0, \frac{1}{k}]}) =\infty$ for any $k\geq 1.$

\endproclaim

\demo{Proof} Assume that each version of an orthogonal process $X$ on $B$ is not a.e. continuous at $t= 0$. Let us fix a countable dense subset $B_1$ in $B$, $0\in B_1$. An assumposition that the restriction $X_{|B_1}$ is a.e. continuous at $t =0$ leads to the existence of a version $\tilde X$ of $X$, a.e. continuous at $t =0$ and defined on the whole of $B$. Thus, for some $\epsilon >0$,
$$||\sup_{t\in B_1\cap [0, \frac{1}{k}]} |X(t)|  || \geq \epsilon$$
for any $k\geq 1$, and
$$||\sup_{t\in B_{1, k}} |X(t)| \ || \geq \frac{\epsilon}{2}$$
for some finite sets $B_{1, k}\subset B_1\cap [0, \frac{1}{k}]$, and $X_{|B_2}$ is not a.e. continuous at $t= 0$ for a countable set
$$B_2 =\bigcup_{k\geq 1} B_{1, k} \cup \{0, 1\}$$
satisfying (13). Theorem 1.7. implies then inequality $V H_B\geq V H_{B_2} = V h_{B_2} =\infty.$

Assume now that $V(H_B 1_{[0, \frac{1}{k}]}) =\infty$ for any $k\geq 1$. Let $B^{ir}$ be a set of 'irregular' points of $B$:
$$t\in B^{ir}\iff t\in B\wedge (t\notin \overline{B\cap (t, \infty)} \vee t\notin \overline{B\cap (-\infty, t)}).$$
Then $\lambda(B^{ir}) =0$ and $V(H_B 1_{(0, 1]\setminus B^{ir}}) =\infty$. Let us observe that for
$$\{0, 1\}\subset B_1\subset B_2\subset \dots \subset B\setminus B^{ir},\quad
\sharp B_l <\infty,$$
with  $\bigcup_{l\geq 1} B_l$ being a dense (countable) subset of $B\setminus B^{ir}$, we have
$$h_{B_l}\nearrow H_B \quad\text{on } (0, 1]\setminus B^{ir}.$$
In particular, $\lim_{l\to\infty} V(h_{B_l} 1_{[0, k]}) =\infty$ for any $k$ and the existence of a subset $\Cal B^s\subset B$ satisfying (13), with $V h_{B^s} =\infty$ is pretty obvious. By Theorem 1.7, there exists an a.e. discontinuous at $0$ process  $X_1$ on $B^s$. Any extension of $\tilde X_1$ to an orthogonal process $X$ on $B$ has then the property that any version of $X$ is a.e. discontinues at $t=0$.
\enddemo

\medskip
\flushpar{\bf 6.7. Proof of Theorem 6.2.} Take first $t\in [0, 1]\setminus \{n 3^{-2^j}; j\geq 0, 0\leq n \leq 3^{2^j}\}$. Then Lemmas 6.3, 6.4, 6.5 reduce the problem to Lemma 6.6.

Assume now that $t = n 3^{-2^j}$ for some $j\geq 0$, $0\leq n\leq 3^{2^j}$. For any process $X$ orthogonal on $B$, the process $X(1-t) - X(1)$ is obviously orthogonal on $1 -B =\{1 -t; t\in B\}$. Thus Lemma 6.6 can be proved, equally well, in the following version.

A'. Any orthogonal process on $B$ is a.e. continuous at $t =1$ if and only if $V(H_B 1_{[\frac{k-1}{k}, 1]}) <\infty$ for some $k\geq 1.$

B'. There exists an orthogonal process on $B$ a.e. discontinuous at $t= 1$ if and only if $V(H_B 1_{[\frac{k-1}{k}, 1]}) =\infty$ for any $k \geq 1.$

We can assume now that $t = n 3^{-2^j}$, $0 < n < 2^j$. For $S_t(s) = [s -t]$, with $[s] =\max (-\infty, s)\cap \Bbb Z$ being the integer part, we have $S_t(t) =0$ and the following equivalences become completely elementary:

Any orthogonal process on $B$ is a.e. continuous at $t$ if and only if any orthogonal process on $S_t(B)$ is a.e. continuous at $0$ and at $1.$

There exists an orthogonal process on $B$ a.e. discontinuous at $t$ if and only if there exists an orthogonal process on $S_t(B)$ a.e. discontinuous at $0$ or a.e. discontinuous at $1$.

Assuming  $[t -\frac{1}{k}, t+ \frac{1}{k}]\in [0, 1]$, we have
$$V(H_B\cap 1_{[t-\frac{1}{k}, t+\frac{1}{k}]}) < \infty$$
if and only if 
$$V(H_{S_t(B)}\cap 1_{[0, \frac{1}{k}]}) <\infty$$
and
$$ V(H_{S_t(B)} \cap 1_{[\frac{k-1}{k}, 1]}) <\infty.$$

Lemma 6.3 is thus proved for any $t\in [0, 1].$
\hfill $\square$

\example{6.8. Example} Any orthogonal process on the Cantor set $C\subset [0, 1]$ has an a.e. continuous version.

Indeed, for the information function $H_C$ we have the estimates
$$||(H_C - k)^+||\leq \sum_{l\geq k} (\frac{2}{3})^k = 3(\frac{2}{3})^k$$
and $\sum_{i\geq 0} ||(H_C -2^i)^+|| <\infty$, hence obviously $V H_C <\infty.$

\endexample

We add some remarks on orthogonal measures. Let $(\Omega, \Cal F, P)$ be a probability space with $\Omega$ a (countable) union of atoms $\Omega = \bigcup_{n\geq 1}\Omega_n$. Then $H(\omega) = -\log_3 P(\Omega_n)$ for $\omega \in\Omega_n$ is, up to a constant factor, the classical information function on $\Omega$. 
Assume for simplicity, that $H\geq 1$.

\definition{6.9. Definition} A function $m: \Cal F\to L_2(\Bbb R)$ is an orthogonal measure if

i)\quad \ $||m(A)||^2 = P(A)$ for $A\in\Cal F;$

ii)\quad $m(\bigcup_{n\geq 1} A_n) =\sum_{n\geq 1} m(A_n)$
in $L_2$-norm, for any  mutually disjoint sets $A_1, A_2,\dots\in\Cal F.$

\enddefinition

As a refolmulation of the  classical results of Moric and Tandori, we have

\proclaim{6.10. Theorem} {\rm A.} The equality $m(\bigcup_{n\geq 1} A_n) =\sum_{n\geq 1} m(A_n)$ is a.e. valid for any orthogonal measure and any disjoint sets $A_1, A_2,\dots \in\Cal F$ if and only if $\sum_{i\geq 0} ||H_i|| <\infty.$

{\rm B.} A sequence $m(B_n)$ is almost everywhere divergent for some orthogonal measure $m$ and some $B_1\supset B_2\supset \dots$ in $\Cal F$ if and only if $\sum_{i\geq 0} ||H_i||=\infty.$

\endproclaim

As before $H_i = H\wedge 2^{i+1} - H\wedge 2^i$, $i\geq 0.$

\demo{Proof} Any sequence $m(B_n)$, $n\geq 1$, for $B_1\supset B_2\supset\dots$, $\bigcap_{n\geq 1} B_n=\emptyset$, can be identified with some tails
$$\sum_{p\geq k(n)} m(\Omega_{\sigma(p)}),\quad n\geq 1,$$
of some rearranged series $\sum_{n\geq 1} m(\Omega_{\sigma(n)})$.

Part $A$ is thus a consequence of Theorem 1.7 B.

Part B is also easy to obtain, using 1.7 B.
\enddemo

\medskip

\flushpar{\bf 6.11. Stationary processes.} For a stationary process (in weak sense) $X: \Bbb R\to L_2$, the investigation of a.e. continuity requires a special approach. It is only suggested by our results on information function. For example
$$\sum_{k\geq 0} (2^k(4||X(2^{-k-1}) - X(0)||^2
- ||X(2^{-k}) - X(0)||^2)^{\frac12} <\infty$$
proves to be a sufficient condition, in some sense the weakest possible. More details will be given in a subsequent paper.

\bye
\proclaim{5.7. Lemma} For any Borel functions $g$, $g_0,\dots, g_i\geq 0$, we have
$$\multline ||V_0(g_0+V_1)\dots (g_i +V_{i+1}) g||\\
\leq ||V_0\dots V_{i+1} g|| +\sum_{0\leq j\leq i} 2^{(j+1)/2}||g_j||.\endmultline$$
\endproclaim

\demo{Proof} To start our backward induction, let us observe that
$$\multline ||V_i(g_i +V_{i+1})g -  V_iV_{i+1} g|| \leq ||(g_i +V_{i+1} g)\wedge 2^i - (V_{i+1} g) \wedge 2^i||\\
+|| \ ||(g_i +V_{i+1} g-2^i)^+||_i - ||(V_{i+1} g-2^i)^+||_i||\leq ||(g_i +V_{i+1}g)\wedge 2^i-V_i V_{i+1} g\wedge 2^i||\\
+ ||(g_i +V_{i+1} g-2^i)^+ - (V_{i+1} g-2^i)^+||\leq ||g^1|| +||g^2||\endmultline$$
for some $g^1, g^2\geq 0$, $g^1+g^2 = g_i$. Thus
$$||V_i (g_i +V_{i+1})g - V_iV_{i+1} g||\leq \sqrt{2} g_i.$$
Then
$$||(g_{i-1}+V_i)(g_i +V_{i+1}g -V_iV_{i+1}g||\leq \sqrt{2}||g_{i-1}|| + \sqrt{2}^2 ||g_i||$$
and so on. Finally, we have (25).
\enddemo

Let us fix a finite set $B$ satisfying (7) and let, for simplicity, $h_B \geq 1.$

\proclaim{5.8. Lemma} Let $\tilde B$ be a generated set and let $h_{\tilde B} < 2^{i+1}$. Then
$$||V_0\dots V_{i+1} \underset{=\!=}\to{h_{\tilde B}}|| \leq ||V_0\dots V_{i+1} h_B||+7.$$

\endproclaim

\demo{Proof} Let us denote $g = h_B$ and
$$g_j = ((\underset{=\!=}\to{h_{\tilde B}} \vee h_B)\wedge 2^{j+1}-2^j)^+ - (h_B \wedge 2^{j+1} -2^j)^+,\quad 0\leq j\leq i.$$
Then
$$(\!={=}\to{h_{\tilde B}} -2^i)^+ =0$$
and
$$||V_0\dots V_{i+1} \underset{=\!=}\to{h_{\tilde B}}||
\leq ||V_0\dots V_{i+1} (\underset{=\!=}\to{h_{\tilde B}} \vee h_B)||\leq ||V_0(g_0+V_1)\dots (g_i +V_{i+1})h_B||.$$
Thus, By Lemma 5.6, we have
$$||V_0\dots V_{i+1} \underset{=}\to{h_{\tilde B}}||\leq ||V_0\dots V_{i+1} h_B||+\sum_{0\leq j\leq i} 2^{(j+1)/2}||g_j||.$$

Let us fix an  interval $(\alpha, \beta]$ satisfying $[\alpha, \beta]\cap B =\{\alpha, \beta\}$ and
$3^{-2^{l+1}}\leq \beta-\alpha <3^{-2^l}$, thus obviously $l+1\leq \log_2 \max_{0 <t \leq 1} h_B(t)\leq \log_2 \max_{0 < t\leq 1} h_{\tilde B}(t) \leq i+1.$

The main idea of the proof is contained in an observation that
$$\lambda((\underset{=}\to{h_{\tilde B}} \geq 2^{j+1})\cap (\alpha, \beta])
\leq \lambda((\underset{=}\to{h_{\tilde B}}\geq 2^{j+2})\cap (\alpha, \beta])\leq 2\cdot 3^{-2^{j+2}}$$
for $l\leq j<i$, and
$$(\underset{=}\to{h_{\tilde B}} < 2^{l+1}) = (\underset{=}\to{h_{\tilde B}} \leq 2^l),$$
$$(\underset{=}\to{h_{\tilde B}} <2^l)\cap (\alpha, \beta]\subset (\underset{=}\to{h_{\tilde B}} \leq h_B)\cap (\alpha, \beta].$$
In other words,
$$g_j 1_{(\alpha, \beta]} ||^2 \leq 2^{2j}\cdot 2\cdot 3^{-2^{j+1}}$$
for $j\geq l+1,$
$$||g_j 1_{(\alpha, \beta]}|| =0 \quad\text{for } \ j\leq l.$$
Estimating the number of mutually disjoint interval $(\alpha, \beta]$, $\beta -\alpha \geq 3^{-2^{l+1}}$, we obtain
$$||g_j||^2 \leq 2^{2j}\cdot 2\cdot 3^{-2^{j+1}}\cdot 3^{2^{l+1}} =2^{2j+1} 3^{-2^j}$$
and
$$\sum_{0\leq j\leq i} 2^{(j+1)/2} ||g_j|| <6.$$
\enddemo

The proof of Theorem 1.7 is now reduced to an implication $V h_B <\infty\implies V\underset{=}\to{h_{\tilde B}} <\infty$. Unexpectively, for arguments being not necesarily triadic functions, like $h_B$, $\underset{=}\to{h_{\tilde B}} $, the comparison of values $V_i h_B$, $V_i \underset{=}\to{h_{\tilde B}} $ of our crucial operation $V_i$ is complicated. To compress calculations, we developed convention 1.4.

\medskip
\flushpar{\bf 5.11. Notation.} For any function $f\geq 1$ we write
$$_{j\downarrow} = f\wedge 2^j$$
and
$$f_j = f_{j+1\downarrow} -f_{j\downarrow}, \quad f_{j\uparrow} = f -f_{j\downarrow},$$
for $j\geq 0$, in paricular
$$f = f_{j\downarrow} + f_j + f_{j+1\uparrow}$$
for any real function $f.$

\proclaim{5.12. Lemma} For any Borel functions $1\leq g\leq h$ on $(0, 1]$ and for
$$a =\bigcup_{\smallmatrix 0\leq n < 3^{2^j}\\ \delta_n^j\cap(h\geq 2^j)\ne\emptyset\endsmallmatrix} \delta_n^j\tag31$$
we have
$$\multline ||(V_j h)_{j-1\uparrow} - (V_j g)_{j-1\uparrow}||\\
 \leq ||h_{j\uparrow} -g_{j\uparrow}|| + ||(h_{j\downarrow} -g_{j\downarrow}) 1_A|| + ||h_{j-1}- g_{j-1}||.\endmultline$$

\endproclaim

\demo{Proof} Demanded calculations are natural:
$$\aligned (V_jh)_{j-1\uparrow} & - (V_j g)_{j-1\uparrow}\\
& = (h_{j-1} - g_{j-1})1_A c + ((V_j h)_{j-1\uparrow} - (V_j g)_{j-1\uparrow}) 1_A\\
& \leq (h_{j-1} - g_{j-1}) 1_A c + (V_j h - V_jg) 1_A\\
&\leq (h_{j-1}- g_{j-1} 1_A c + ||h_{j\uparrow} - g_j\uparrow}||_j + (h_{j\downarrow} -g_{j\downarrow})1_A\\
& \leq (h_{j-1}-g_{j-1}) + ||h_{j\uparrow} -g_{j\uparrow}|| + (h_{j\downarrow} -g_{j\downarrow}) 1_A.\endaligned$$
\enddemo

\proclaim{5.13. Lemma} For any set 4B$ satisfying (13) and for $0\leq j <i$ let us denote
$$g = V_{j+1}\dots V_i h_B,$$
$$h = V_{j+1} \dots V_i (h_B\vee \underset{=}\to{h_{\tilde B}} )$$
and let the set $A$ be giveb by (31). Then
$$||h_{j-1} - g_{j-1}|| \leq 2^j\cdot 3^{-2^j},\tag32$$
$$||(h_{j\downarrow} -g_{j\downarrow}) 1_A||\leq 2^{j+1} 3^{-2^{j-2}}.\tag33$$
\endproclaim

\demo{Proof} The inequality (32) is an immediate consequence of our Definition 5.10 giving  of $\underset{=}\to{h_{\tilde B}} $. Let us fix an interval $(\alpha, \beta]$ satisfying $[\alpha, \beta]\cap B = \{\alpha, \beta\}$. Observe that
$$(\alpha, \beta]\cap (\underset{=}\to{h_{\tilde B}}  >2^{j-1}) = (\alpha, \beta]\cap (h_{\tilde B} > 2^{j+1})\subset \delta_k^{j+1}\cup\delta_l^{j+1}$$
for some $0\leq k\leq l < 3^{2^{j+1}}$, and, because of the structure of $V_{j+1},\dots, V_i$,
$$(\alpha, \beta]\cap (h > 2^{j-1}) \subset \delta_k^{j+1}\cup \delta_l^{j+1}.$$
Thus $\lambda((\alpha, \beta]\cap [h > 2^{j-1})) < 2\cdot 3^{-2^{j+1}}.$

On the other hand
$$\multline \sharp \{(\alpha, \beta];\quad [\alpha, \beta]\cap B =\{\alpha, \beta\}, (g < 2^j)\cap (\alpha, \beta]\ne \emptyset\}\\
\leq \sharp \{(\(\alpha, \beta]; \quad [\alpha, \beta]\cap B =\{\alpha, \beta\}, (h_B <2^j)\cap (\alpha, \beta]\ne\emptyset\}\leq 1/3^{-2^j},\endmultline$$
as $(h_B)_{j+1\downarrow} \leq g_{j+1\downarrow}$. Inequality (32) is given by a natural estimation
$$||h_{j-1}-g_{j-1}||\leq 2^{j-1}\bigcup_{\smallmatrix (\alpha, \beta]\\ [\alpha, \beta] \cap B =\{\alpha, \beta\}\endsmallmatrix} \lambda((\alpha, \beta]\cap (g <2^j)\cap (h >2^{j-1})).$$

To prove (33) let us fix $\delta_m^j\subset A$, $0\leq m < 3^{2^j}$. Supposing $\delta_m^j\cap B =\emptyset$ we have $\delta_m^j\cap \tilde B =\emptyset$ and $\delta_m^j\subset A^c$. A careful analysis is needed only when
\newpage
$$k =\min \{n; \delta_n^{j+1}\cap\delta_m^j\cap B\ne \emptyset\},$$
$$l= \max \{n; \delta_n^{j+1}\cap \delta_m^j\cap B\ne \emptyset\}$$
are defined. The case $k= l$ is obviously possible. Then 
$$h_B >2^j \quad\text{on } \bigcup_{k <n<l} \delta_n^{j+1}$$
and, by the structure of $V_{j+1},\dots, V_i$,
$$h\geq g >2^j\quad\text{on } \bigcup_{k <n < l} \delta_n^{j+1}.$$
Moreover $h, g$ are constant on the set $\delta_m^j\cap \bigcup_{n <k}\delta_n^{j+1}$ as well as on the set $\delta_m^j\cap \bigcup_{n >l} \delta_n^{j+1}$,
$$h \leq 2^{j-1}\quad\text{on } \delta_m^j\cap \bigcup_{n <k \ \text{or } n >l} \delta_n^{j+1}.$$
We say that $m$ is special if
$$(g < 2^{j-2})\cap \delta_m^j \cap \bigcup_{n <k \ \text{or } n >l} \delta_n^{j+1} \ne \emptyset.$$
Then, by natural estimation of number of large intervals $[\alpha, \beta)$, $[\alpha, \beta]\cap B =\{\alpha, \beta\}$, $\beta -\alpha >3^{-2^{j-1}}$, we have
$$\multline||h_{j\downarrow} - g_{j\downarrow}||^2 \leq 2\cdot 2^{2j} \cdot 3^{-2^{j+1}} \cdot \sharp \{0,\dots, 3^{2^j}-1\}\\
+ 2^{2j}\cdot 3^{-2^j}\cdot\sharp\{m; n \ \text{is special}\}\leq 2\cdot 2^{2j}(3^{-2^j} +3^{-2^j}\cdot 3^{2^{j-1}}),\endmultline$$
thus (33).\hfill $\square$
\enddemo

\medskip

\flushpar{\bf 5.13. Proof of Theorem 1.7}. For any set $B$ satisfying (13) with $h_B\leq K < \infty$ we have $V h_B\leq K$ for any subset $B_1\subset B$, satisfying (13). For suitable large $i \geq 1$, Lemmas 5.11, 5.12 give inequalities
$$\aligned  ||V_0\dots & V_i \underset{==}\to{h_{\tilde B_1}} -1||\\
& \leq ||V_0\dots V_i(\underset{=}\to{h_{\tilde B_1}} \vee h_{B_1}) -1||\\
&\leq ||V_0\dots V_i h_{B_1} -1|| +\sum_{j\geq 0} (2^j\cdot 3^{-2^j} + 2^{j+1}\cdot 3^{-2^{j-2}})\\
& \leq K +\sum_{j\geq 0} (2^j\cdot 3^{-2^j}+2^{j+1}\cdot 3^{-2^{j-2}}) =: L,\endaligned$$
and, by Lemma 5.10,
$$||V_0\dots V_i \underset{=\!=}\to{h_{\tilde B_1}} || \leq 14(L +1).$$
Then $V h_{\tilde B} < 14(L +1)$ by an obvious application of Lebesgue theorem on monotonic convergence and Lemmas 5.7, 4.6 gives a.e. convergence of any orthogonal process on $B$.

Together with 4.9, it completes the proof of 1.7.

$$The proof of Theorem 1.7 is now reduced to an implication $V h_B <\infty\implies V\underset{=}\to{h_{\tilde B}} <\infty$. Unexpectively, for arguments being not necesarily triadic functions, like $h_B$, $\underset{=}\to{h_{\tilde B}} $, the comparison of values $V_i h_B$, $V_i \underset{=}\to{h_{\tilde B}} $ of our crucial operation $V_i$ is complicated. To compress calculations, we developed convention 1.4.

\medskip
\flushpar{\bf 5.11. Notation.} For any function $f\geq 1$ we write
$$_{j\downarrow} = f\wedge 2^j$$
and
$$f_j = f_{j+1\downarrow} -f_{j\downarrow}, \quad f_{j\uparrow} = f -f_{j\downarrow},$$
for $j\geq 0$, in paricular
$$f = f_{j\downarrow} + f_j + f_{j+1\uparrow}$$
for any real function $f.$

\proclaim{5.12. Lemma} For any Borel functions $1\leq g\leq h$ on $(0, 1]$ and for
$$a =\bigcup_{\smallmatrix 0\leq n < 3^{2^j}\\ \delta_n^j\cap(h\geq 2^j)\ne\emptyset\endsmallmatrix} \delta_n^j\tag31$$
we have
$$\multline ||(V_j h)_{j-1\uparrow} - (V_j g)_{j-1\uparrow}||\\
 \leq ||h_{j\uparrow} -g_{j\uparrow}|| + ||(h_{j\downarrow} -g_{j\downarrow}) 1_A|| + ||h_{j-1}- g_{j-1}||.\endmultline$$

\endproclaim

\demo{Proof} Demanded calculations are natural:
$$\aligned (V_jh)_{j-1\uparrow} & - (V_j g)_{j-1\uparrow}\\
& = (h_{j-1} - g_{j-1})1_A c + ((V_j h)_{j-1\uparrow} - (V_j g)_{j-1\uparrow}) 1_A\\
& \leq (h_{j-1} - g_{j-1}) 1_A c + (V_j h - V_jg) 1_A\\
&\leq (h_{j-1}- g_{j-1} 1_A c + ||h_{j\uparrow} - g_j\uparrow}||_j + (h_{j\downarrow} -g_{j\downarrow})1_A\\
& \leq (h_{j-1}-g_{j-1}) + ||h_{j\uparrow} -g_{j\uparrow}|| + (h_{j\downarrow} -g_{j\downarrow}) 1_A.\endaligned$$
\enddemo

\proclaim{5.13. Lemma} For any set 4B$ satisfying (13) and for $0\leq j <i$ let us denote
$$g = V_{j+1}\dots V_i h_B,$$
$$h = V_{j+1} \dots V_i (h_B\vee \underset{=}\to{h_{\tilde B}} )$$
and let the set $A$ be giveb by (31). Then
$$||h_{j-1} - g_{j-1}|| \leq 2^j\cdot 3^{-2^j},\tag32$$
$$||(h_{j\downarrow} -g_{j\downarrow}) 1_A||\leq 2^{j+1} 3^{-2^{j-2}}.\tag33$$
\endproclaim

\demo{Proof} The inequality (32) is an immediate consequence of our Definition 5.10 giving  of $\underset{=}\to{h_{\tilde B}} $. Let us fix an interval $(\alpha, \beta]$ satisfying $[\alpha, \beta]\cap B = \{\alpha, \beta\}$. Observe that
$$(\alpha, \beta]\cap (\underset{=}\to{h_{\tilde B}}  >2^{j-1}) = (\alpha, \beta]\cap (h_{\tilde B} > 2^{j+1})\subset \delta_k^{j+1}\cup\delta_l^{j+1}$$
for some $0\leq k\leq l < 3^{2^{j+1}}$, and, because of the structure of $V_{j+1},\dots, V_i$,
$$(\alpha, \beta]\cap (h > 2^{j-1}) \subset \delta_k^{j+1}\cup \delta_l^{j+1}.$$
Thus $\lambda((\alpha, \beta]\cap [h > 2^{j-1})) < 2\cdot 3^{-2^{j+1}}.$

On the other hand
$$\multline \sharp \{(\alpha, \beta];\quad [\alpha, \beta]\cap B =\{\alpha, \beta\}, (g < 2^j)\cap (\alpha, \beta]\ne \emptyset\}\\
\leq \sharp \{(\(\alpha, \beta]; \quad [\alpha, \beta]\cap B =\{\alpha, \beta\}, (h_B <2^j)\cap (\alpha, \beta]\ne\emptyset\}\leq 1/3^{-2^j},\endmultline$$
as $(h_B)_{j+1\downarrow} \leq g_{j+1\downarrow}$. Inequality (32) is given by a natural estimation
$$||h_{j-1}-g_{j-1}||\leq 2^{j-1}\bigcup_{\smallmatrix (\alpha, \beta]\\ [\alpha, \beta] \cap B =\{\alpha, \beta\}\endsmallmatrix} \lambda((\alpha, \beta]\cap (g <2^j)\cap (h >2^{j-1})).$$

To prove (33) let us fix $\delta_m^j\subset A$, $0\leq m < 3^{2^j}$. Supposing $\delta_m^j\cap B =\emptyset$ we have $\delta_m^j\cap \tilde B =\emptyset$ and $\delta_m^j\subset A^c$. A careful analysis is needed only when
$$k =\min \{n; \delta_n^{j+1}\cap\delta_m^j\cap B\ne \emptyset\},$$
$$l= \max \{n; \delta_n^{j+1}\cap \delta_m^j\cap B\ne \emptyset\}$$
are defined. The case $k= l$ is obviously possible. Then 
$$h_B >2^j \quad\text{on } \bigcup_{k <n<l} \delta_n^{j+1}$$
and, by the structure of $V_{j+1},\dots, V_i$,
$$h\geq g >2^j\quad\text{on } \bigcup_{k <n < l} \delta_n^{j+1}.$$
Moreover $h, g$ are constant on the set $\delta_m^j\cap \bigcup_{n <k}\delta_n^{j+1}$ as well as on the set $\delta_m^j\cap \bigcup_{n >l} \delta_n^{j+1}$,
$$h \leq 2^{j-1}\quad\text{on } \delta_m^j\cap \bigcup_{n <k \ \text{or } n >l} \delta_n^{j+1}.$$
We say that $m$ is special if
$$(g < 2^{j-2})\cap \delta_m^j \cap \bigcup_{n <k \ \text{or } n >l} \delta_n^{j+1} \ne \emptyset.$$
Then, by natural estimation of number of large intervals $[\alpha, \beta)$, $[\alpha, \beta]\cap B =\{\alpha, \beta\}$, $\beta -\alpha >3^{-2^{j-1}}$, we have
$$\multline||h_{j\downarrow} - g_{j\downarrow}||^2 \leq 2\cdot 2^{2j} \cdot 3^{-2^{j+1}} \cdot \sharp \{0,\dots, 3^{2^j}-1\}\\
+ 2^{2j}\cdot 3^{-2^j}\cdot\sharp\{m; n \ \text{is special}\}\leq 2\cdot 2^{2j}(3^{-2^j} +3^{-2^j}\cdot 3^{2^{j-1}}),\endmultline$$
thus (33).\hfill $\square$
\enddemo

\medskip

\flushpar{\bf 5.13. Proof of Theorem 1.7}. For any set $B$ satisfying (13) with $h_B\leq K < \infty$ we have $V h_B\leq K$ for any subset $B_1\subset B$, satisfying (13). For suitable large $i \geq 1$, Lemmas 5.11, 5.12 give inequalities
$$\aligned  ||V_0\dots & V_i \underset{=}\to{h_{\tilde B_1}} -1||\\
& \leq ||V_0\dots V_i(\underset{=}\to{h_{\tilde B_1}} \vee h_{B_1}) -1||\\
&\leq ||V_0\dots V_i h_{B_1} -1|| +\sum_{j\geq 0} (2^j\cdot 3^{-2^j} + 2^{j+1}\cdot 3^{-2^{j-2}})\\
& \leq K +\sum_{j\geq 0} (2^j\cdot 3^{-2^j}+2^{j+1}\cdot 3^{-2^{j-2}}) =: L,\endaligned$$
and, by Lemma 5.10,
$$||V_0\dots V_i \underset{=}\to{h_{\tilde B_1}} || \leq 14(L +1).$$
Then $V h_{\tilde B} < 14(L +1)$ by an obvious application of Lebesgue theorem on monotonic convergence and Lemmas 5.7, 4.6 gives a.e. convergence of any orthogonal process on $B$.

Together with 4.9, it completes the proof of 1.7.

\bye

\proclaim{Lemma 3.20} For any function $h_i\in \Cal T_i$, $i\geq 5$, there exists an operator $U$ of type $i$ such that
$$(U h_i -2^i)^+ \geq (h_i -2^i)^+ - f_i -g_i,$$
for some positive functions $f_i, g_i$,
$$f_i + g_i\leq (h_i -2^i)^+,$$
$$||f||_i \leq 2^{-i},$$
$$g_i \leq 2^{-i}(h_i -2^{i})^+.$$
\endproclaim

Unfortunately, to obtain an implication (10) for $V {h_B} \geq C y$, $y >1$, we cannot omit some more special properties of triadic functions (passages 3.22, 3.23, 3.24). Some simplification are given by the following convention.

For any bounded Borel function $f$ on $(0, 1]$ and any operation $U: L_2 (0, 1]\to L_2(0, 1]$ we use $(f + U)$ as a symbol for new operation  $(f +U) h = U h+f$.

Let us observe that for any triadic function $h$ we have
$$ V_j h, V_j(-g +h) \quad\text{are triadic},\quad
h\wedge 2^j = V_j h\wedge 2^j = (V_j (-g +h))\wedge 2^j,$$
for any Borel function $g$, $0\leq g\leq (h -2^j)$, and for $j\geq 0.$

\proclaim{3.22. Lemma} For any triadic function $h$ and any positive Borel functions $f_1,\dots, f_i$, $g_1\dots, g_1$ on $(0, 1]$ satisfying
$$f_i +g_i \leq (V_{i+1} h-2^i)^+\tag12$$
and
$$f_j+g_j\leq (V_{j+1}(-f_{j+1}-g_{j+1}+V_{j+2})\dots (-f_i-g_i+V_{i+1})h-2^j)^+,\quad 1\leq  k_j <i,\tag13$$
we have
$$||V_1(-f_1-g_1+V_2)\dots(-f_i-g_i +V_{i+1})h-V_1(-g_1+V_2)\dots (-g_i+V_{i+1}h||\leq ||f_1||+\dots+||f_i||.$$
\endproclaim

\demo{Proof} Nothing by subadditivity  of $L_2$-norm is used. For any triadic functions $h_1, h_2$ with $h_1\wedge 2^j =h_2 \wedge 2^j$, we have
$$V_j h_1-V_j h_2 =||(h_1-2^j)^+ ||_j - ||(h_2-2^j)^+||_j\leq ||h_1-h_2||_j,$$
and
$$||V_j h_1 -V_j h_2|| \leq ||h_1-h_2||.$$

Let us denote
$$d_j = (-g_j +V_{j+1})\dots (-g_i +V_{i+1})h,$$
$$e_j = (-f_j -g_j +V_{j+1})\dots (-f_i -g_i +V_{j+1})h,$$
then $e_j\wedge 2^j =d_j\wedge 2^j$, and $h_1:=V_j d_j$, $h_2 := V_jd_j$ are triadic function with $h_1\wedge 2^j= h_2\wedge 2^j$.

We use backward induction. Observe at first that
$$||V_i d_i-V_ie_i|| = || \ ||d_i-2^i||_i - ||e_i -2^i||_i||
\leq || \ || d_i-e_i||_i|| = ||d_i-e_i|| = ||f_i||.$$
Then supposing that
$$||V_{j+1} d_{j+1} -V_{j+1} e_{j+1}||\leq ||f_{j+1}||+\dots + ||f_i||$$
for some $1\leq j <i$, by analogical calculations, we obtain
$$\multline ||V_jd_j -V_j e_j||\leq ||V_j(-f_j -g_j +V_{j+1})d_{j+1}\\
- V_j(-g_j+V_{j+1}) e_{j+1}||\leq ||f_j|| +\dots + ||f_i||.\endmultline$$

Finally $||V_1d_1 - V_1e_1||\leq ||f_1|| +\dots + ||f_i||$ and the lemma is proved.

\enddemo

\proclaim{3.23. Lemma } For a triadic function $h$, Borel functions $g_j$ and for $0 <\epsilon_j <1$, $1\leq j\leq l$, if we have
$$0\leq g_i \leq \epsilon_i(V_{i+1} h-2^i)^+\tag14$$
and $$0\leq g_j \leq \epsilon_j(V_{j+1} (-g_{j+1} +V_{j+1})\dots (-g_iV_{i+1}) h-2^j)^+, \quad 1\leq j <i,\tag15$$
then
$$||V_0V_1(-g_1+V_2)\dots (-g_i+V_{i+1}h||\geq (1-\epsilon_1)\dots (1-\epsilon_i) ||V_0\dots V_{i+1} h||.\tag14$$
\endproclaim

\demo{Proof}  By definition of $V_i$, and by monotonicity of $|| \cdot||_i, V_i$,
$$\multline||(V_i (-g_i +V_{i+1})h-2^{i-1})^+||_{i-1}=||((h\wedge 2^i)-2^{i-1})^+ +||(-g_i+V_{i+1}h-2^i)^+||_i||_{i-1}\\
\geq ||((h\wedge 2^i)-2^{i-1})^+ +(1-\epsilon_i)||(V_{i+1}h-2^i)^+||_i||_{i-1}\\
\geq ||(1-\epsilon_i)[((h\wedge 2^i)-2^{i-1})^+ + ||(V_{i+1} h-2^i)^+||_i]||_{i-1}=\epsilon_i||(V_iV_{i+1}h-2^{i-1})^+||_{i+1}.\endmultline$$
Supposing that
$$\multline||V_j(-g_j+V_{j+1})\dots (-g_i+V_{i+1})h-2^{j-1}||_{j-1}\\
 \geq (1-\epsilon_j)\dots(1-\epsilon_i) ||(V_j\dots V_{i+1} h-2^{j-1})^+||_{j-1}\endmultline$$
for some $1<j\leq i$, we have analogously
$$\multline ||(V_{j-1} (-g_{j-1} +V_j)\dots (-g_i +V_{i+1})h-2^{j-2}||_{j-2}\\
 = ||(h\wedge 2^{j-1} -2^{j-2})^+
+ ||((-g_{j-1} +V_j)\dots (-g_i +V_{i+1}) h-2^{j-1})^+||_{j-1}
||_{j-2}\\
 \geq ||(h\wedge 2^{j-1} -2^{j-2})^+
+ (1-\epsilon_{j-1})||(V_j (-g_j +V_{j+1})\dots (-g_i +V_{i+1})h -2^{j-1})^+||_{j-1}||_{j-2}\\
 \geq ||(h\wedge 2^{j-1}-2^{j-2})^+
+ (1-\epsilon_{j-1})\dots (1-\epsilon_i) ||(V_j\dots V_{i+1} -2^{j-1})^+||_{j-1}||_{j-2}\\
\geq (1-\epsilon_{j-1})\dots (1-\epsilon_i)||(V_{j-1}\dots V_{i+1} h-2^{j-2})^+||_{j-2}.\endmultline$$
\bye
By backward induction, we have in particular
$$||V_1 (-g_1+V_2)\dots (-g_i+ V_{i+1}) h-1||_0 \leq (1-\epsilon_1)\dots (1-\epsilon_i)||(V_1\dots V_{i+1} h-1)^+||_0,$$
thus (14).

\enddemo

\proclaim{3.24. Corollary} For any triadic function $h$ and any positive Borel functions $f_1,\dots, f_i, g_1,\dots, g_i$ on $(0, 1]$ satisfying (12), (13) and (14), (15) we have
$$\multline ||V_0V_1(-f_1-g_1 +V_2)\dots (-f_i-g_i +V_{i+1}) h|| \\
\geq (1-\epsilon_1)\dots (1-\epsilon_i) ||V_0\dots V_{i+1} h|| - ||f_1|| -\dots -||f_i||.\endmultline$$
\endproclaim

\demo{Proof} By Lemma 3.22 we have also
$$\multline||V_0 V_1(-f_1-g_1+V_2)\dots (-f_i -g_i +V_{i+1})h||\\
\geq ||V_0V_1 (-g_1 +V_2)\dots (-g_i +V_{i+1}) h|| - ||f_1|| -\dots -||f_i||\endmultline$$
and Lemma 3.23 can be used.

\enddemo

\medskip

\flushpar{\bf 3.25. The proof of Basic Lemma 3.3}. Being precise, we  show the implication (10)
$$||V h_B|| \geq C y\quad\text{implies } \ [0, 1]\in \Cal S^B(\frac12,  y)$$
for any finite triadic set $B$ and $y >1$. 

Let us assume that $h_B \leq i+1$. By 3.18,
$$[0, 1]\in\Cal S^B(\frac12, ||(V_8 U_8\dots V_i U_i \underline{h_B} -2^8)^+||)$$
for any operators $U_j$ of type $j$, $8 \leq j\leq i$, cf. 3.14, 3.15. Then, by 2.20, using back word induction, one can obtain operators $U_8,\dots, U_i$ in such a way, that
$$(U_j(V_{j+1} U_{j+1}\dots V_iU_i) \underline{h_B} -2^j)^+ \geq (V_{j+1}U_{j+1}\dots V_iU_i \underline{h_B} -2^j)^+-f_j-g_j,$$
for some positive functions $f_j$, $g_j$,
$$f_j+g_j \leq (V_{j+1}U_{j+1}\dots V_iU_i \underline{h_B}-2^j)^+,$$
$$||f_j||\leq 2^{-j},$$
$$g_j\leq 2^{-j}(V_{j+1} U_{j+1}\dots V_i U_i \underline{h_B}-2^j)^+,$$
for all indices $8\leq j\leq i$. For a triadic set $B$ and for $h_B\leq 2^{i+1}$, it is obvious that $V_{j+1}\underline{h_B} = \underline{h_B}$ and Corollary 3.24 can be used, with $\epsilon_1 =\dots =\epsilon_7 =0$, $f_1 =g_1=\dots =f_7 =g_7 =0$, giving
$$\multline ||V_0\dots V_7(V_8U_8)\dots (V_iU_i)\underline{h_B}||\\
\geq (1-2^{-8})\dots (1-2^{-i}) ||V_0\dots V_i\underline{h_B}|| -2^{-8} -\dots -2^{-i}.\endmultline$$
By 3.21, we have
$$||V_0\dots V_i \underline{h_B}||\geq \frac{1}{3} ||V_0\dots V_i h_B|| = \frac{1}{3} V h_B$$
and, by definition of our basic operations $V_j$ (cf. 1.6),
$$\multline ||V_0\dots V_7 (V_8U_8) \dots (V_i U_i) \underline{h_B}||\\
\leq ||1_{(0, 1]}|| + \sum_{0\leq j\leq 7} ||2^j \cdot 1_{(0, 1]}|| + ||(V_8U_8\dots V_i U_i \underline{h_B} -2^8)^+||.\endmultline$$
Thus any $C\geq  \frac{3}{1-2^{-7}} +2^8+2^{-7}$ can be used in (10).

 of discontinuous process on a giveb set $B$ satisfying (7), with $V h_B=\infty$, can be deduced by rather formal and simple way.

First of all the existence of a discontinuous process on $B$ is equivalent to the existence of a  discontinuous process on a triadic, "generated set" $\tilde B$ (Lemma 4.6).

Obviously $h_{\tilde B} \ge q H_B$ and the condition $V h_{\tilde B} =\infty$ is satisfied. Thus, for any given $C >0$, intervals of the form
$$\delta_s = [3^{-2^{i(s+1)}}, 3^{-2^{i(s)}}],\quad i(1) < i(2) <\dots, \ 3^{-2^{i(s)}}\in\tilde B,\tag16$$
can be taken in such a way that
$$Vh_{\tilde B\cap \delta_s} >C.$$
Monotonicity of operations $V_j$ and a Lebesgue limit theorem are natural tools.

By Tshebysthev inequality more specific form of our basic Lemma 3.3 can be obtained

\proclaim{4.2. Corollary} There exists a constant $C >0$ such that inequality 
$$V(h_{\tilde B} 1_\delta) \geq C$$
for
$$\delta = (3^{-2^i}, 3^{-2^j}],\quad 0\leq j< i,\quad 3^{-2^i}, 3^{-2^j}\in{\tilde B}$$
implies
$$\lambda(\max_{t\in{\tilde B}\cap \delta} X(t) >1) >\frac{1}{4}$$
for some orthogonal process $X: {\tilde B}\cap \delta \to L_2[0, 1]$ (i.e. $\Bbb E X(t) =0$, $\Bbb E(X(t) - X(s))^2 = (t-s)$, for $s, t\in{\tilde B}\cap \delta$, $X(3^{-2^i}) =0$).
\endproclaim

The proof we give letter. Now, a suitable orthogonal process $X$ on $\tilde B$ with independent increments $X(t_k) - X(s_k)$, $k=1,2,\dots$ for $s_k, t_k\in \tilde B\cap \delta_k$ with $\delta_k$ given by (16) is a.e. discontinuous. It is given by a standard use of Borel-Cantelli Lemma.

\medskip

\flushpar{\bf 4.3. Outlook of the proof of continuity of a process}. Assume now that $V h_B <\infty$ for some set $B$ satisfying (7). It is enough to show that each orthogonal process $X$ on $\tilde B$ is continuous a.e.

We can not omitted some tedious estimations to show that $V h_B <\infty$ implies $V h_{\tilde B} <\infty$. But then the result is almost immediate. We construct, by backward induction some majorants $M_{ij} \in L_2$,
$$\sup_{t\in\tilde B\cap [3^{-2^i}, 3^{-2^j}]} |X(t)| \leq M_{ij}, \quad 0\leq j<i,$$
with $||M_{ij}||^2 \leq K_j <\infty$, and $K_j \searrow 0$.
A crucial role is played by classical Plancheree idea of diadic partitions of a gives sequence of vectors. It gives

\proclaim{4.4. Lemma (see \ )} For any othogonal vectors $Y_1,\dots Y$ in $L_2$, we have
$$||M||^2 \leq k^2 \sum_{1\leq n\leq N} ||X||^2$$
for
$$M =\max_{1\leq n\leq N} |Y_1 +\dots + Y_N|,$$
$$k \geq \log_2 N.$$
\endproclaim

We pass to a detail construction of a discontinuous orthogonal process $X$ on $B$ for $V h_B =\infty$. Let $\rho(t, B) = \inf_{s\in B} |s -t|$. We start with the following gemetrical observation

\bigskip
\flushpar{\bf 4.7. Proof of Corollary 4.2}. For $B_\delta =\widehat{\tilde B\cap \delta}$ and for any orthogonal process $Y$, $Y: B_\delta\to L_2[0, 1)$, $Y(0) =0$, we have
$$\multline ||\max_{t\in\tilde B\cap\delta} (Y(t) - Y(3^{-2^i}))||\\
\geq ||\max_{t\in B_\delta} Y(t)|| - (\sum_{t\in B_\delta} \rho(t, \tilde B\cap\delta))\geq ||\max_{t\in B_\delta} Y(t)|| - \sqrt{3},\endmultline$$
by Lemma 4.4.

Thus
$$\lambda(\max_{t\in\tilde B\cap\delta} (X(t) - Y(3^{-2^i})) \leq 1) <\frac{3}{4}$$
if only
$$\lambda(\max_{t\in B_\delta} Y(t) <1 + 2\sqrt{3}) <\frac12.$$
As $h_{B_\delta}\geq h_{\tilde B} 1_\delta$, our Lemma 3.3 can be used.

\bigskip

Elementary properties of conditional $L_2$-norm, as
$$||h_k+ h_0||_i \underset{k\to\infty}\to{\longrightarrow} \infty \quad\text{implies } ||h_k||_i
\underset{k\to\infty}\to{\longrightarrow} \infty$$
for any bounded function $h_0$, leads to analogical properties of operations $V_i$ and in consequence
$$V(h\wedge 2^i) \underset{i\to\infty}\to{\longrightarrow} \infty \quad\text{for } V h =\infty,$$
$$V(h- h_0) =\infty \quad \text{for } V h =\infty, \ 0< h_0 \leq h\wedge 2^j.$$
Thus intervals $\delta_s$ of the form (16), satisfying $V h_{\tilde B\cap\delta_s} >C$, can be obtain. Let $X_s(t)$, $s=1,2,\dots,$ be processes given by 4.2 with $\delta =\delta_s$.

We can not omitted some tedious estimations to show that $V h_B <\infty$ implies $V h_{\tilde B} <\infty$. But then the result is almost immediate. We construct, by backward induction some majorants $M_{ij} \in L_2$,
$$\sup_{t\in\tilde B\cap [3^{-2^i}, 3^{-2^j}]} |X(t)| \leq M_{ij}, \quad 0\leq j<i,$$
with $||M_{ij}||^2 \leq K_j <\infty$, and $K_j \searrow 0$.
A crucial role is played by classical Plancheree idea of diadic partitions of a gives sequence of vectors. It gives

We pass to a detail construction of a discontinuous orthogonal process $X$ on $B$ for $V h_B =\infty$. Let $\rho(t, B) = \inf_{s\in B} |s -t|$. We start with the following gemetrical observation

\bye

As previously, the investigation of a finite triadic set $B$ is a main step:

\proclaim{5.1. Lemma} For any finite traidic set $B$, with $h_B\leq 2^{i+1}$, and any orthogonal process $X(t)$, $t\in B$, we have
$$||M^0||\leq ||2 \log_2 3(V_0,\dots, V_i \underline{h_b} -1)||\tag21$$
for $M^0 = \max_{t\in B} |X(t)|.$
Moreover,
$$M^j||\leq ||2\log_2 3(V_j\dots V_i\underline{h_B} -2^j 1_{[0, 3^{-2^j}]}||\tag22$$
for $3^{-2^j}\in B$ and
$$M^j =\max_{t\in B, t\leq 3^{-2^j}} |X(t)|.$$
\endproclaim

\demo{Proof} Let us denote
$$M_m^j =\max_{t\in B\cap\delta_m^j} |X(t) -X(m 3^{-2^j})|$$
if $(m 3^{-2^j}, (m+1) 3^{-2^j})\cap B\ne \emptyset$, and $M_m^j =0$ otherwise. As before
$$\delta_m^j =[m 3^{-2^j}, (m+1) 3^{-2^j}]$$
and we take $0\leq j\leq i$, $0\leq m < 3^{2^j}.$

For a moment we use additional notations $\bar 1=1$,
$$\bar a = 2^{j+1}\quad\text{for } 2^j < a \leq 2^{j+1}, \ j\geq 0,$$
$$\bar V_j h =h\wedge 2^{j+1} + ||(h -2^{j+1})^+||_j.$$
To describe a step of backward induction, assume that, for some triadic function $h_j$,
$$M_n^{j+1}||\leq \log_2 3||(\bar h_j 1_{\delta_n^{j+1}}-2^{j+1})^+||^2$$
for any $0\leq n < 3^{2^{j+1}}$. Then always
$$M_m^j|| \leq \log_2 3||(\bar V_j(\bar h_j 1_{\delta_m^j}) -2^j)^+||$$
for $0\leq m < 3^{2^i}$. Indeed, for $\delta_m^j \subset (h_j\leq 2^j)$, we have $\delta_m^j\cap B=\emptyset$ and $M_m^j =0$. For $\delta_m^j\notin (h_j\leq 2^j)$, we assume that $X(t)$ is extended to an orthogonal process with $t\in B\cup (\delta_m^j\cap \{\delta_m^j\cap \{n 3^{-2^{j+1}}; 0\leq n \leq 3^{2^j}\})$.  Then
$$\multline ||M_m^j|| \leq\max_{t\in \delta_m^j\cap \{n 3^{-2^{j+1}}; 0\leq n\leq 3^{2^{j+1}}\}} |X(t) - X(m 3^{-2^j})| \ ||\\
+ (\sum_{\underset{\delta_n^{j+1}\subset\delta_m^j}\to{0\leq n \leq 3^{2^{2^{j+1}}}}} ||M_n^{j+1}||^2)^\frac12\leq \log_2 3||2^j 1_{\delta_m^j}||\\
+(\sum_{\underset{\delta_n^{j+1}\subset \delta_m^j}\to{0\leq n\leq 2^{j+1}}} ||(\bar h 1_{\delta_n^{j+1}} -2^{j+1})^+||^2)^\frac12= ||(\bar V_j\cdot \bar h_j\cdot 1_{\delta_m^j} -2^j)^+||.\endmultline$$

The multiplier $\log_2 3$ appear because of Lemma 4.4. Moreover $\overline{\bar V_j\cdot\bar h_j} = \bar V_j \bar h_j$, and by backward induction we have
$$||M^0|| = \log_2 3||\bar V_0\dots \bar V_i \bar h_B -1||.$$

Observe now that
$$\bar V_j \bar h = 2 V_j \frac12 \bar h \leq 2 V_j \underline{h}$$
for any triadic function $h$, and $\overline{2 V_{\underline{h}}} = 2 V\underline{h}.$

Thus inequality (21) is valid. The proof of (22) is almost identical.

\enddemo

\bye

We pass to a detail construction of a discontinuous orthogonal process $X$ on $B$ for $V h_B =\infty$. Let $\rho(t, B) = \inf_{s\in B} |s -t|$. We start with the following gemetrical observation

\bye